\documentclass[authoryear, 11pt]{imsart}
\RequirePackage[OT1]{fontenc}
\RequirePackage{amsthm,amsmath,natbib}
\startlocaldefs
\numberwithin{equation}{section}
\theoremstyle{plain}

\endlocaldefs
\usepackage{psfrag, epsfig, graphicx}
\usepackage{amsfonts,amsmath,latexsym,amssymb}
\usepackage[active]{srcltx}
\newtheorem{lemma}{Lemma}
\newtheorem{theorem}{Theorem}
\newtheorem{proposition}{Proposition}

\newtheorem{assumption}{Assumption}
\newtheorem{remark}{Remark}
\newtheorem{definition}{Definition}

\newtheorem{assumption*}{Assumption}

\renewcommand{\kappa}{\varkappa}

\newcommand{\rd}{{\rm d}}
\newcommand{\e}{\varepsilon}

\newcommand{\cA}{{\cal A}}
\newcommand{\cB}{{\cal B}}

\newcommand{\cE}{{\cal E}}
\newcommand{\cF}{{\cal F}}
\newcommand{\cG}{{\cal G}}
\newcommand{\cH}{{\cal H}}
\newcommand{\cI}{{\cal I}}

\newcommand{\cK}{{\cal K}}

\newcommand{\cP}{{\cal P}}
\newcommand{\cQ}{{\cal Q}}
\newcommand{\cR}{{\cal R}}

\newcommand{\cU}{{\cal U}}

\newcommand{\cZ}{{\cal Z}}

\newcommand{\Bg}{\boldsymbol{g}}

\newcommand{\bp}{\boldsymbol{p}}

\newcommand{\bq}{\boldsymbol{q}}
\newcommand{\BQ}{\boldsymbol{Q}}

\newcommand{\bld}{\boldsymbol{\delta}}

\newcommand{\bE}{\mathbb E}

\newcommand{\bH}{\mathbb H}

\newcommand{\bK}{{\mathbb K}}
\newcommand{\bL}{{\mathbb L}}

\newcommand{\bN}{{\mathbb N}}

\newcommand{\bP}{{\mathbb P}}
\newcommand{\bQ}{{\mathbb Q}}
\newcommand{\bR}{{\mathbb R}}

\newcommand{\bV}{{\mathbb V}}

\newcommand{\mF}{\mathfrak{F}}

\newcommand{\mI}{\mathfrak{I}}

\newcommand{\mh}{\mathfrak{h}}

\newcommand{\mm}{\mathfrak{m}}
\newcommand{\mn}{\mathfrak{n}}

\newcommand{\mH}{\mathfrak{H}}

\newcommand{\ma}{\mathfrak{a}}
\newcommand{\mb}{\mathfrak{b}}

\newcommand{\mz}{\mathfrak{z}}
\newcommand{\mQ}{\mathfrak{Q}}
\newcommand{\rp}{\mathrm{p}}

\newcommand{\epr}{\hfill\hbox{\hskip 4pt
                \vrule width 5pt height 6pt depth 1.5pt}\vspace{0.5cm}\par}

\newcommand{\blh}{\boldsymbol{h}}

\newcommand{\mfi}{\mathbf{i}}

\newcommand{\bn}{\mathbf{n}_\mathbf{i}}

\begin{document}
\begin{frontmatter}
\title{Structural adaptation in the density model.}
\runtitle{Structural adaptation}
\begin{aug}
\author[t1]
{\fnms{O.V.} \snm{Lepski}
\ead[label=e1]{oleg.lepski@univ-amu.fr}}
\author[t1]{\fnms{G.} \snm{Rebelles}
\ead[label=e2]{rebelles.gilles@neuf.fr}}
\thankstext{t1}{This work has been carried out in the framework of the Labex Archim\`ede (ANR-11-LABX-0033) and of the A*MIDEX project (ANR-11-IDEX-0001-02), funded by the "Investissements d'Avenir" French Government program managed by the French National Research Agency (ANR).}
\runauthor{O.V. Lepski and G. Rebelles}

\affiliation{Aix--Marseille Universit\'e, CNRS, Centrale Marseille, I2M UMR 7373}

\address{Institut de Math\'ematique de Marseille\\
Aix-Marseille  Universit\'e   \\
 39, rue F. Joliot-Curie \\
13453 Marseille, France\\
\printead{e1}\\
\printead{e2}}
\end{aug}
\begin{abstract}

This paper deals with non-parametric  density estimation on $\bR^2$ from i.i.d observations.
It is assumed that after \textsf{unknown} rotation of the coordinate system the coordinates
of the observations are independent random variables whose densities belong to a H\"older class with
\textsf{unknown} parameters.  The minimax and adaptive minimax theories for this structural statistical
model are developed.

\end{abstract}


\end{frontmatter}

\section{Introduction}

Let $\xi\in\bR^2$ be a random vector having the density $g$ w.r.t the Lebesgue measure. We will assume that the coordinates of $\xi$ are independent
and let $X\in\bR^2$ be the random vector obtained from the relation

\vskip0.0cm

\centerline{$
X=M\xi,\quad M\in\mQ,
$}

\vskip0.1cm

\noindent where $\mQ$ is the set of all  rotational  $2\times2$-matrices.

\noindent Let we observe $n\in\bN^*$ independent copies of $X$ that is $X^{(n)}=(X_1,\ldots,X_n)$. We want to estimate the density of $X$ denoted by $f$ at a given point $x\in\bR^2$ using the observations $X^{(n)}$.
 By estimator, we mean any $X^{(n)}$-measurable map $\widetilde{f}:\bR^n\to \bR$. The accuracy of an estimator $\hat{f}$
is measured by the pointwise risk
\vskip0.1cm
\centerline{$
 \cR^{(p)}_n[\widetilde{f}, f]:=\Big(\bE_{f} \big|\widetilde{f}-f(x)\big|^p\Big)^{1/p},\;p\in [1,\infty).
$}
\vskip0.1cm
\noindent Here $\bE_{f}$ denotes the expectation with respect to the probability measure
$\bP_{f}$ of the observations $X^{(n)}$.

Let $\cQ\subseteq\mQ$ be fixed and let $\cG(\beta,L)$ denote the following set of functions.

\begin{definition}
\label{def:cG} We say that  $g:\bR^2\to \bR$ belongs to $\cG(\beta,L)$ if
\begin{enumerate}

\item[($\mathbf{i}$)] $g(\cdot,\cdot)=g_1(\cdot)g_2(\cdot)$ and $g_1,g_2:\bR\to\bR_+$ are symmetric probability densities;

\smallskip

\item[($\mathbf{ii}$)] $g_1,g_2$  belong to the H\"older class $\bH(\beta,L)$, $\beta>0,L>0$, on $\bR$.

\end{enumerate}
\end{definition}
\noindent  For the reader's convenience the formal definition of $\bH(\beta,L)$ is postponed to the end of this section.  Here we only mention that $\beta$ is referred to the smoothness of the underlying function while $L$ is the Lipschitz constant.

\noindent For any $\beta>0,L>0$ introduce the following set of probability densities.

\vskip0.2cm

\centerline{$
\cF(\beta,L,\cQ)=\left\{f:\bR^2\to\bR_+:\; f(\bullet)=g\big(M^T\bullet\big),\; g\in\cG(\beta,L),\; M\in\cQ\right\}.
$}

\vskip0.2cm

\noindent In the present paper we will study the minimax and minimax adaptive estimation of the density $f$ over  the collection of functional classes
$\cF(\beta,L,\cQ)$.
To illustrate the interesting feature of the  problem at hand let us consider the simplest situation. Assume that the set $\cQ$ consists a single element $\BQ$. In this case we can first obtain new observation sequence $\xi_1=\BQ^{T}X_1,\ldots,\xi_n=\BQ^{T}X_n$. Noting that the density of $\xi_1$ is $g_1g_2$
we estimate next separately $g_1$ and $g_2$ from the sequence of the first and second coordinates of $\xi_1,\ldots,\xi_n$ respectively. In particular one can use the kernel estimation method with properly chosen bandwidth. It will lead to the estimators $\widehat{g}_1$ and  $\widehat{g}_2$. Since $g_1,g_2\in\bH(\beta,L)$ the pointwise minimax accuracy (minimax rate of convergence) of each marginal density  will be proportional to $n^{-\frac{\beta}{2\beta+1}}$. Therefore, the minimax pointwise accuracy  in estimating of $g$ provided by the estimator $\widehat{g}(x)=\widehat{g}_1(x_1)\widehat{g}_2(x_2)$ is proportional to $n^{-\frac{\beta}{2\beta+1}}$ as well.
The estimator for $f(x)=g\big(\BQ^Tx\big)$ is then given by $\widehat{f}_{\BQ}(x)=\widehat{g}(\BQ^Tx)$.

All saying above can be summarized as follows.
\begin{theorem}
\label{th1}
Let $\beta>0$, $L>0$ and $\mathbf{\BQ\in\mQ}$ be fixed. Then, for any $x\in\bR^2$ there exists an estimator $\widehat{f}_{\BQ}(x)$  such that $\forall p\geq 1$

\vskip0.15cm

\centerline{$
\displaystyle{\sup_{\cF(\beta,L,\{\BQ\})} \cR^{(p)}_n[\widehat{f}_{\BQ}(x), f]\lesssim n^{-\frac{\beta}{2\beta+1}}}.
$}

\vskip0.15cm

\noindent Moreover (here and later $\inf$ is taken over all possible estimators) $\forall p\geq 1$

\vskip0.15cm

\centerline{$
\displaystyle{\inf_{\widetilde{f}}\sup_{\cF(\beta,L,\{\BQ\})} \cR^{(p)}_n[\widetilde{f}, f]\gtrsim n^{-\frac{\beta}{2\beta+1}}}.
$}

\end{theorem}
\vskip-0.2cm
\noindent The proof of this theorem is straightforward. Moreover its first assertion  follows from the results obtained in Proposition \ref{prop2} presented in Section \ref{sec:proofs}.

The assertions of Theorem \ref{th1} show that the structural assumption $f(\bullet)=g\big(M^T\bullet\big)$ leads to the essential improvement of the accuracy of estimation. Indeed, it is easily seen that  $\cF(\beta,L,\mQ)\subset\bH(\vec{\beta},\vec{L})$, where $\bH(\vec{\beta},\vec{L})$ is the isotropic H\"older class  on $\bR^2$ with $\vec{\beta}=(\beta,\beta)$ and $\vec{L}=(L^2,L^2)$. Recall that the minimax pointwise accuracy on this class is given by
$n^{-\frac{\beta}{2\beta+2}}$ which is much larger than the univariate rate $n^{-\frac{\beta}{2\beta+1}}$ available under the structural assumption discussed above.

The first problem which we address is the following: do the statements of Theorem \ref{th1} remain valid if the cardinality of $\cQ$ is larger than $1$? We remark that  the matrix $M$ describing the law of observation is unknown in this case. Therefore, we are talking about the adaptation to unknown rotation of coordinate system (structural adaptation). We will see that the answer on aforementioned question depends heavily on the "massiveness" of the set $\cQ$.
In particular, Theorem \ref{th1} is not valid if $\cQ=\mQ$. On the other hand if $\cQ$ is a finite set whose elements satisfy some separation condition and their number is independent of $n$ the assertions of Theorem \ref{th1} hold.

\noindent The second problem studied in the paper is the minimax adaptive estimation with respect to the parameter $(\beta,L)$. Let $\cQ\subseteq\mQ$ be fixed and let

\vskip0.15cm

\centerline{$
\displaystyle{\varphi_n(\beta,L)=\inf_{\widetilde{f}}\sup_{\cF(\beta,L,\cQ)} \cR^{(p)}_n[\widetilde{f}, f],\quad \beta>0, \;L>0}.
$}

\vskip0.1cm

\noindent
Our objective is to answer on the following question: does there exist an estimator $f^*$ such that

\vskip0.15cm

\centerline{$
\displaystyle{\limsup_{n\to\infty}\varphi^{-1}_n(\beta,L)\sup_{\cF(\beta,L,\cQ)} \cR^{(p)}_n[f^*, f]<\infty,\quad\forall\beta>0, \;L>0?}.
$}

\vskip0.1cm

\noindent
We will prove that the answer is positive if $\cQ$ is a net in $\mQ$ satisfying some separation condition and $\beta\in (0,\mb]$, where $\mb>0$ is an arbitrary but a priori chosen number.

\paragraph{Historical notes}
There is a vast literature dealing with minimax and minimax adaptive density estimation. The interested reader can find  very detailed
overview on this topic in \cite{lepski15}.
As it was saying above,  we will follow the modeling strategy which consists in imposing
additional structural assumptions on the function to be estimated. This approach
was pioneered by \cite{stone} who discussed the trade-off between flexibility and
dimensionality of nonparametric models and formulated the heuristic dimensionality
reduction principle. Standard examples of structural nonparametric models are
single-index,  additive, projection pursuit or multi-index model,  composite functions structure etc.
The minimax and minimax adaptive results in these models (mostly in the nonparametric regression context) were obtain in \cite{huber}, \cite{chen}
\cite{golubev}, \cite{hris}, \cite{horowitz}, \cite{juditsky}, \cite{GL09}, \cite{o+n} among many others.
However, when one is talking about the multivariate density estimation there are not so many articles  where minimax and  minimax adaptive results were obtained. The problems and models similar to those considered in  the present paper were studied in   \cite{samarov}, \cite{amato}, \cite{lepski13a}, \cite{rebelles1},
\cite{rebelles2}.
We would like especially to mention the paper \cite{samarov1} where $d$-dimensional variant of our model was considered. Some problems in this article  have been studied
 under pointwise risk and  we will provide a detailed comparison of them and our results after Theorem \ref{th3}.

\paragraph{Definitions, assumptions and notations} For any $\cQ\in\mQ$ and any function $f\in\cF(\beta,L,\cQ)$ we denote by $\BQ_f\in\cQ$ and $\Bg_f\in\cG(\beta,L)$ the quantities obtained from the relation

\vskip0.0cm

\centerline{$
\displaystyle{f(\bullet)=g_f\big(\BQ_f\bullet\big)}.
$}

\vskip0.1cm

\noindent
Obviously this representation  is not unique and later on we consider an arbitrary couple
 $(\BQ_f,\Bg_f)$ for which the latter relation holds.

Furthermore $\|\cdot\|_\infty$ will be used for the supremum norm on $\bR$, the integer part of  $a>0$ will be denoted by $\lfloor a\rfloor$
  and
any $Q\in\mQ$ will be presented as

\vskip0.1cm

\centerline{$
\displaystyle{Q=(q,q_\perp)=\left(\begin{array}{cc}
q_1  & -q_2
\\
q_2  & q_1
\end{array}
\right)}.
$}
\begin{definition}
\label{de:holder}
Let $\beta=r+\alpha, r\in\bN$, $0<\alpha\leq 1$ and $L>0$ be given. We say that $w:\bR\to\bR$ belongs to the H\"older class $\bH(\beta,L)$ if it is $r$-times continuously differentiable, $\|w^{(j)}\|_\infty\leq L$ for any $j=0,\ldots,r$ and

\vskip0.2cm

\centerline{$
\displaystyle{\|w^{(r)}(\cdot+\mz)-w^{(r)}(\cdot)\|_\infty\leq L|\mz|^\alpha,\quad\forall \mz\in\bR}.
$}

\vskip0.0cm

\end{definition}
\noindent
For given $\mb\geq 1$ we denote by $\bK_\mb$ the set of functions
 $\cK:\bR\to \bR$ satisfying the following assumption.
\begin{assumption}
\label{ass:kernel}
$\cK\in\bL_1(\bR)\cap\bL_\infty (\bR)$, $\int_\bR \cK(u)\rd u=1$ and

\vskip0.2cm

\centerline{$
\int_{\bR} \cK(u)u^{j}\rd u=0, j=1,\ldots,2\lfloor\mb \rfloor,\quad \int_{\bR}|\cK(t)||t|^{2\mb}\rd t<\infty.
$}

\vskip0.0cm

\end{assumption}

\noindent With any $\cK\in\bK_\mb$ we associate the following quantity:

\vskip0.2cm

\centerline{
$C(\cK,\mb,\mathbf{s})=\sup_{b\leq\mb}\int_{\bR^2} \big|\cK(t_1)\cK(t_2)\big|\big[\mathbf{s}\big(t_1^2+t_2^2\big)^{b}+1\big]^{2}\rd t_1\rd t_2,\; \mathbf{s}>0.
$}

\vskip0.2cm

\noindent For any $D,Q\in\mQ$  we will write

\vskip0.2cm

\centerline{$
p_1:=p_1(D,Q)=q^Td_\perp,\quad\; p_2:=p_2(D,Q)=q^Td
$}

\vskip0.2cm

\noindent and set
$
\varrho(D,Q):=\min\big[|p_1(Q,D)|,|p_2(D,Q)|\big].
$

\noindent For given $\delta>0$ we denote by $\bQ_\delta$ the set of all subsets of $\mQ$ consisting of $\delta$-distinguishable points with respect to $\varrho$.
Recall that $Q_1,Q_2$ are called $\delta$-distinguishable  with respect to $\varrho$ if
$
\varrho(Q_1,Q_2)\geq \delta.
$
For any $\cQ_\delta\in\bQ_\delta$ let

\vskip0.15cm

\centerline{$
\mn(\cQ_\delta)=\ln{\big(\text{card}(\cQ_\delta)\big)}.
$}

\begin{remark}
\label{rem1}

 We note that $p_1(D,Q)=-p_1(Q,D)$, $p_2(D,Q)=p_2(Q,D)$ and $p_1(Q,Q)=0$. Additionally it can be easily checked that
 \vskip0.15cm

\centerline{$
 \varrho(Q_1,Q_3)\leq 2\sqrt{2}\big[\varrho(Q_1,Q_2)+\varrho(Q_2,Q_3)\big],\quad\forall Q_1,Q_2,Q_3\in\mQ.
 $}

 \vskip0.15cm

\noindent Hence, we assert that $\varrho$ is a $2\sqrt{2}$-pseudo-inframetrics   on $\mQ$.

\end{remark}
\vskip-0.2cm
\noindent If $\mQ_\delta\in\bQ_\delta$ is the maximal $\delta$-net of $\mQ$ in $\varrho$ then   $\mn_\delta:=\mn(\mQ_\delta)$ is called $\bld$-\textsf{capacity} of $\mQ$.
Recall that the $\delta$-capacity (as well as the $\delta$-entropy) is used for   classifying  compact  metric  sets  according  to  their  massivity.

\noindent From now on $\delta\in(0,1)$ (possibly dependent on $n$) is assumed to be fixed
and the number of observations $n\geq 3$.


\section{Main results}
\label{sec:main}
In this section we develop  the minimax and adaptive minimax theories over  collection of functional classes $\cF(\beta,L,\cQ_\delta), \cQ_\delta\in\bQ_\delta$.

\subsection{Lower bounds} We start with presenting two lower bound results.
\begin{theorem}
\label{th:lb}
For any $\beta_1>0$, $\beta_2>0$ and $p\geq 1$ there exists $\mathbf{c}_1>0$ such that for any  $\BQ\in\mQ$
and $L>0$
\begin{equation}
\label{eq:lb-adaptive}
\liminf_{n\to\infty}\;\inf_{\widetilde{f}}\sup_{\beta\in\{\beta_1,\beta_2\}}\big(L^{\frac{2}{2\beta+1}}\ln(n)/n\big)^{-\frac{\beta}{2\beta+1}}\sup_{\cF(\beta,L,\{\BQ\})} \cR^{(p)}_n[\widetilde{f}, f]\geq \mathbf{c}_1.
\end{equation}
For any $\beta>0$,  $L>0$ and $p\geq 1$ there exists $\mathbf{c}_2>0$ such that for an arbitrary sequence $\delta_n>0$ satisfying $\delta_n\geq \big(\ln(n)\big)^{\frac{2\beta+2}{2\beta+1}}\big(L^{-2}/n\big)^{\frac{1}{2\beta+1}}$ and any $\cQ_{\delta_n}\in\bQ_{\delta_n}$ one has

\vskip0.15cm

\centerline{$
\displaystyle{\liminf_{n\to\infty}\big(\mn(\cQ_{\delta_n})/n\big)^{-\frac{\beta}{2\beta+1}}\;\inf_{\widetilde{f}}\;\sup_{\cF(\beta,L,\cQ_\delta)} \cR^{(p)}_n[\widetilde{f}, f]\geq \mathbf{c}_2}.
$}

\end{theorem}
\vskip-0.3cm

Some remarks are in order.

 $\mathbf{1^0}.$ The first assertion of the theorem is quite standard and its proof  will be omitted.  The fact that  $\BQ$ is known  reduces the considered problem to adaptive poinwise estimation over collection of H\"older classes under independent hypothesis (it suffices to consider new observation sequence $\BQ^{T}X_1,\ldots,\BQ^{T}X_n$). Then (\ref{eq:lb-adaptive}) follows, in particular, from  the lower bound result obtained in \cite{rebelles1}.
As usual, see for instance \cite{rebelles1},  there is $\ln(n)$-price to pay for adaptation. That means that minimax result given in Theorem \ref{th1} differs from whose in (\ref{eq:lb-adaptive}) by  $\ln(n)$-factor.

 $\mathbf{2^0}.$ The proof of the second assertion is much more involved. If  $\delta=constant$ then the factor $\mn(\cQ_{\delta})$ can be viewed as the \textsf{price to pay for structural adaptation} (with respect to unknown rotation $\BQ_f\in\cQ_{\delta}$). However this is a \textsf{constant factor}, the asymptotics of minimax risk with respect to $n$ remains the same and coincides with whose  in Theorem \ref{th1}. The situation changes completely if $\delta=\delta_n\to 0$, $n\to\infty$. Indeed, if $\mn(\cQ_{\delta_n})\to\infty$ the minimax rate found in Theorem \ref{th1} is no more achievable and $\mn(\cQ_{\delta_n})$ is the minimal price to pay for structural adaptation over $\cQ_{\delta_n}$.
 It is not difficult to see that for any $\cQ_{\delta_n}\in\bQ_{\delta_n}$
  \begin{equation}
 \label{eq:bound-n_delta}
 \mn(\cQ_{\delta_n})\leq\mn(\mQ_{\delta_n})\asymp\;|\ln(\delta_n)|,\;n\to\infty.
 \end{equation}
 It yields in particular that if $\delta_n\sim n^{-a}$ for some $a>0$, then the minimal price to pay for structural adaptation on $\mQ_{\delta_n}$ is proportional to $\ln(n)$.

\subsection{Pointwise selection rules}
\label{sec:selection-rule}

 Our estimation procedures are based on the original selection rule from the family of kernel-type estimators.
 One of them called  \textsf{adaptive selection rule} is inspired by general approach discussed in \cite{GL12} but the procedure is completely new.

\smallskip

\noindent {\it Family of estimators.} For any $\cK$ satisfying Assumption \ref{ass:kernel} and $h>0$ denote  $\cK_h(\cdot)=h^{-1}\cK(\cdot/h)$.
 For any $D\in\mQ$ and $x\in\bR^2$
 introduce the estimator

\vskip0.1cm

\centerline{$
\widetilde{f}_{h,D}(x)=\Big[n^{-1}\sum_{k=1}^n\cK_{h}\big(d^T(X_k-x)\big)\Big]\Big[n^{-1}\sum_{k=1}^n\cK_{h}\big(d^T_\perp(X_k-x)\big)\Big].
$}

\vskip0.1cm

\noindent Set $\cH=\big\{e^{-k},\; k=0,1,\ldots, \lfloor\ln(n)\rfloor\big\}$ and let  $\cQ_\delta\in\bQ_\delta$ be given. Introduce the following estimator's family:

\vskip0.15cm

\centerline{$
\mF(\cQ_\delta,\cH)=\big\{\widetilde{f}_{h,D}(x),\;D\in\cQ_\delta,\; h\in\cH\big\}.
$}

\vskip0.15cm

\noindent It is worth noting that if $\cQ_\delta=\{D\}$  the estimator $\widetilde{f}_{h,D}(x)$ is exactly  the estimator $\widehat{g}(D^Tx)$, $\widehat{g}(x)=\widehat{g}_1(x_1)\widehat{g}_2(x_2)$ introduced in the discussion preceded Theorem \ref{th1}.

Below we propose two different  data-driven selection rules from this collection. The first one, called below \textsf{adaptive selection rule}, will be used in the situation when the parameters $\beta,L$ are unknown,  $\cQ_{\delta_n}$ is an arbitrary element  of $\bQ_{\delta_n}$ with any $\delta_n>0$ satisfying
\begin{equation}
\label{eq:delta}
\delta_n\geq (\ln(n)/n)^{\frac{1}{4\mb+2}}.
\end{equation}
Here $\mb\geq 1$ is an arbitrary but a priori chosen number and $\beta\in (0,\mb]$.

\noindent The second one, called \textsf{minimax selection rule}, will be applied when $\beta,L$ are known. The  interesting case here is $\delta=constant$ for example
$card(\cQ_\delta)=2$. Another intriguing case is $\delta=\delta_n$ such that $\mn(\cQ_{\delta_n})=o(\ln(n))$, $n\to\infty$.

\paragraph*{Auxiliary estimator}  Set  $K_h(t)=\cK_h(t_1)\cK_h(t_2),\; t\in\bR^2,\; h>0$,

\vskip0.15cm

\centerline{$
 \Gamma =\left(
\begin{array}{cccc}
 1&  0
\\*[0mm]
0  & -1
\end{array}\right),\quad\;\Omega=\left(
\begin{array}{cccc}
 0&  1
\\*[3mm]
 1  & 0
\end{array}\right),
$}

\vskip0.15cm

\noindent
and define for any $D,Q\in\mQ$ what we will call the \textsf{auxiliary estimator}

\vskip0.2cm

\centerline{$\displaystyle{
\overline{f}_{h,(D,Q)}(x)=\frac{1}{n(n-1)}\sum_{k,l=1,k\neq l}^nK_{h}\big(p_1\Omega\Gamma X_k+p_2X_l-\Omega\Gamma QD\Omega x\big)}.
$}

\vskip0.15cm

\noindent
Remark that $\overline{f}_{h,(D,Q)}$ is a $U$-statistics of the order 2 if $\varrho(D,Q)\neq 0$. Put also

\vskip0.1cm

\centerline{$
\displaystyle{\widetilde{f}_{h,(D,Q)}(x)=\left\{
\begin{array}{cc}
\overline{f}_{h,(D,Q)}(x),\quad& D\neq Q;
\\*[2mm]
\widetilde{f}_{h,Q}(x),\quad& D=Q.
\end{array}
\right.}
$}

\subsubsection{Adaptive selection rule}
\label{sec:selection-rule-adaptive} Let $\mathbf{A}>0$ be a constant given in section \ref{th2:proof}.

\noindent Set $\;\mH=\left\{h\in\cH :\;1/\ln(\ln(n))\geq h\geq \big[\ln(n)\big]^2/n\right\}$ and

\vskip0.1cm

\centerline{$
\displaystyle{\widehat{\cU}_n=\sup_{\eta\in\mH}\sup_{D\in\cQ_{\delta_n}}\sup_{b\in\{d,d_\perp\}}\bigg\{1\vee\bigg[ n^{-1}\sum_{k=1}^n\big|\cK_\eta\big(b^T(X_k-x)\big)\big|\bigg]^2\bigg\}}.
$}

\vskip0.1cm

\noindent
Introduce for any $Q\in\cQ_\delta$ and any $h\in\mH$

\vskip0.1cm

\centerline{$
\displaystyle{R_{n}(Q,h)=\sup_{\stackrel{\eta,\eta^\prime\in\mH:}{\eta^\prime\leq\eta\leq h}}\;\sup_{D\in\cQ_{\delta_n} }\bigg[\big|\widetilde{f}_{\eta,(D,Q)}(x)-\widetilde{f}_{\eta^\prime,D}(x)\big|-\mathbf{A}\widehat{\cU}_n
\bigg(\frac{\ln(n)}{n\eta^\prime}\bigg)^{1/2}
\;\bigg]_+}
$}

\vskip0.0cm

\noindent
and define

\vskip0.1cm

\centerline{$
\displaystyle{\big(\widehat{h,}\widehat{Q}\big)=\arg\min_{Q\in\cQ_{\delta_n}, h\in\mH}\Big[R_{n}(Q,h)+\mathbf{A}\widehat{\cU}_n\sqrt{\ln(n)/nh}\;\Big]}
$}

\vskip0.1cm

\noindent
The suggested estimator is then $\widehat{f}=\widetilde{f}_{\widehat{h}, \widehat{Q}}(x)$.

\begin{theorem}
\label{th2} Let $p\geq 1$, $\mb\geq 1$ and  $\cK\in\bK_\mb$  be fixed.
Then for any $\beta_1,\beta_2\in (0,\mb]$, $L>0$, any  $\delta_n$ satisfying (\ref{eq:delta}) and
 any $\cQ_{\delta_n}\in\bQ_{\delta_n}$ one has

 \vskip0.2cm

\centerline{$
\displaystyle{\limsup_{n\to\infty}\sup_{\beta\in\{\beta_1,\beta_2\}}\big(L^{2/\beta}\ln(n)/n\big)^{-\frac{\beta}{2\beta+1}}\sup_{\cF(\beta,L,\cQ_{\delta_n})} \cR^{(p)}_n\big[\widehat{f}, f\big]<\infty}.
$}
\end{theorem}

\noindent We conclude that the estimator $\widehat{f}$ provides the optimal (in view of the first assertion of Theorem \ref{th:lb}) accuracy of estimation
simultaneously over the collection of
functional classes $\cF(\beta,L,\cQ_{\delta_n})$.

\subsection{Minimax selection rule}
\label{sec:selection-rule-minimax}
As it has been already mentioned the construction of the minimax estimator is much more complicated. In particular  it requires non-trivial splitting of the observation sequence in order to get desirable theoretical results. However the implementation of our minimax procedure for reasonable sample size does not require such splitting, see remark after Theorem \ref{th3}.

Let $\beta>0$, $L>0$, $\delta_n>0$, $\cQ_{\delta_n}\in\bQ_{\delta_n}$ be given.  Introduce the following notations. Set $\ell_\mathbf{0}=\ln(n)$
and let for any $\mathbf{i}\in\bN^*$
\begin{gather*}
\ell_{\mathbf{i}}=\ln{(\ell_{\mfi-1})},\quad \omega_\mfi=\ell_{\mathbf{i}}\vee 4+\mn(\cQ_{\delta_n});
\\
 \mathbf{i}^*=\min\big\{\mathbf{i}\in\bN^*:\; \omega_{\mathbf{i}}=4+\mn(\cQ_{\delta_n})\big\}.
\end{gather*}
Set also  for any $\mfi=1,\ldots\mfi^*-1$

\vskip0.1cm

\centerline{$
\bn=\lfloor n\ell_\mfi^{-1}\rfloor,\quad N_{\mfi}=\lfloor n/4\rfloor+\sum_{\mathbf{j}=1}^{\mfi}\mathbf{n}_{\mathbf{j}}.
$}

\begin{remark}
\label{rem2}
By definition $\ell_{\mfi^*-1}\geq 4$ that obviously implies that $\ell_{\mfi^*-j}\geq 4^{j}$ for any $j=1,\ldots,\mfi^*-1$. Hence

\vskip0.1cm

\centerline{$
N_{\mfi^*-1}\leq \lfloor n/4\rfloor+ n\sum_{\mfi=1}^{\mfi^*-1}\ell^{-1}_\mfi< 3n/4.
$}

\end{remark}
\vskip-0.2cm
\noindent In view of the latter remark introduce the following splitting of the observation sequence. For any $\mfi=1,\ldots\mfi^*-1$ set

\vskip0.2cm

\centerline{$
X^{(\mfi)}=X_{N_{\mfi-1}+1},\ldots,X_{N_\mfi},\;  X^{(\mfi^*)}=X_{N_{\mfi^*-1}+1},\ldots, X_n,\; X^{(\mathbf{0})}=X_{1},\ldots,X_{\lfloor n/4\rfloor}.
$}

\vskip0.2cm

\noindent We remark that $X^{(\mfi)}, \mfi=0,\ldots\mfi^*,$ are mutually independent and  later on all objects measurable with respect to $X^{(\mfi)}, \mfi=1,\ldots\mfi^*$   will be marked by "$^{(\mfi)}$".

\noindent Put  $\mathbf{n}_{\mathbf{i}^*}=n-N_{\mfi^*-1}$ and for any $\mfi=1,\ldots\mfi^*$ introduce $h_{\mathbf{i}}=(L^{-4}\omega_{\mathbf{i}}/\bn)^{\frac{1}{2\beta+1}}$,
\begin{gather*}
R^{(\mfi)}_{n}(Q)=\sup_{D\in\cQ_{\delta_n}}\Big[\big|\widetilde{f}^{(\mfi)}_{h_\mfi,(D,Q)}(x)-\widetilde{f}^{(\mfi)}_{h_\mfi,D}(x)\big|-
\mathbf{B} L^2h_{\mathbf{i}}^{\beta}\;\Big]_+;
\\*[1mm]
\widehat{Q}^{(\mfi)}=\arg\min_{Q\in\cQ_\delta}R^{(\mfi)}_{n}(Q),\quad \widehat{f}^{(\mfi)}(x)=\widetilde{f}^{(\mfi)}_{h_\mfi, \widehat{Q}^{(\mfi)}}(x).
\end{gather*}
\vskip-0.3cm
\noindent where  $\mathbf{B}>0$ is a constant given in section \ref{th3:proof}.

Let  $\widehat{f}(x)$ be the estimator from Theorem \ref{th2} corresponding to the choice
$\mb=\beta$ and constructed from $X^{(\mathbf{0})}$.
Define for any $\mfi=1,\ldots\mfi^*$

\vskip0.1cm

\centerline{$
\displaystyle{\breve{f}^{(\mathbf{0})}(x)=\widehat{f}(x),\qquad\breve{f}^{(\mfi)}(x)=\left\{
\begin{array}{ll}
\breve{f}^{(\mfi-1)}(x),\quad& R^{(\mfi)}_{n}\big(\widehat{Q}^{(\mfi)}\big)\neq 0;
\\*[2mm]
\widehat{f}^{(\mfi)}(x),\quad& R^{(\mfi)}_{n}\big(\widehat{Q}^{(\mfi)}\big)=0.
\end{array}
\right.}
$}

\vskip0.15cm

\noindent
The suggested estimator is then $\breve{f}=\breve{f}^{(\mfi^*)}(x)$.

\begin{theorem}
\label{th3}
Let $p\geq 1$, $\beta>0$, $L>0$,  $\delta_n$ satisfying (\ref{eq:delta}) with $\mb=\beta$, $\cQ_{\delta_n}\in\bQ_{\delta_n}$ and  $\cK\in\bK_\beta$  be fixed.
Then

\vskip0.1cm

\centerline{$
\displaystyle{ \limsup_{n\to\infty}\big(L^{\frac{2}{\beta}}\mn(\cQ_{\delta_n})/n)\big)^{-\frac{\beta}{2\beta+1}}\sup_{\cF(\beta,L,\cQ_{\delta_n})} \cR^{(p)}_n\big[\breve{f}, f\big]<\infty}.
$}

\end{theorem}
\vskip-0.3cm

$\mathbf{1^0}.$ To the best of our knowledge the construction led to the estimator $\breve{f}$ has no analogue in the existing literature on the minimax and minimax adaptive estimation. Although formally
$\mathbf{i}^*\to\infty, n\to\infty$, but $\mathbf{i}^*=1$  for $n=10^{100}$ and for any $\cQ_\delta\in\mQ_\delta$ such that $\mn(\cQ_\delta)\geq 6$.
It worth noting that if $\mathbf{i}^*$ is independent of $n$ (what is the case for many sequences $\delta_n\to 0$) the splitting of data is not needed anymore.
The estimator construction remains the same but the estimators $\widehat{f}^{(\mathbf{i})}, \mathbf{i}=1,\ldots, \mathbf{i}^*$ are built from the whole data set.  The proof of the minimax optimality of this procedure is the simple
modification of the proof of Theorem \ref{th3} and is left to an interested reader.

$\mathbf{2^0}.$ Comparing the results presented  in the second assertion of Theorem \ref{th:lb} and in Theorem \ref{th3} we conclude that $\breve{f}$ is minimax optimal
on $\cF(\beta,L,\cQ_{\delta_n})$. In particular if $\delta>0$ is independent of $n$ there is only a constant factor to be paid   for the  adaptation w.r.t unknown rotation. On the other hand if $\delta_n\asymp n^{-a}, a\leq [4\beta+2)]^{-1}$ and $\mn(\cQ_{\delta_n})\asymp\ln(n)$ for instance $\cQ_{\delta_n}=\mQ_{\delta_n}$, cf. (\ref{eq:delta}), the adaptive estimator from Theorem \ref{th2} with $\mb=\beta$ is minimax optimal as well. We remark that this estimator does not require any splitting of the observations.

$\mathbf{3^0}.$ In \cite{samarov1} the authors studied the same observation model but in an arbitrary dimension $d\geq 2$.  Their estimation procedure is based on the completely different principles. First, they estimate the unknown rotation matrix and then plug-in it to the estimator $\widehat{f}_{\mh,\bullet}$. It is worth noting that the estimation of the rotation (the important problem itself) requires very restrictive assumptions.
In particular the authors assumes that $\beta>5$ if $d=2$ and that the observations possess finite absolute moment of order $4$. We impose none of these assumptions.
Although \cite{samarov1} it is  assumed that the marginal densities belong to  H\"older class the obtained rate of convergence is not uniform one. In particular the authors established the rate which is the same as in Theorem \ref{th1} which, in view of the lower bound of Theorem \ref{th:lb}, is possible if and only if the number of rotations is finite and independent of $n$. We think that the use of plug-in approach in structural models is either too restrictive or not optimal. It seems that the technique of structural adaptation is much more adequate for such kind of problems.

\section{Proofs of Theorems \ref{th:lb}--\ref{th3}}
\label{sec:proofs}
Recall that we will proof only the second assertion of Theorem \ref{th:lb}.

\subsection{Proof of the second assertion of Theorem \ref{th:lb}}
\label{sec-proofs}
To simplicity of  notation we will prove the theorem for $x=0$. The transition to the general case does not bring any additional  difficulty.

 $\mathbf{1^0}.\;$ Let $\mathbf{n}(y)=\frac{1}{\sqrt{2\pi\sigma^2}}e^{-\frac{y^2}{2\sigma}}$, where $\sigma^2>0$ is chosen in order to guarantee $\mathbf{n}(\cdot)\in\bH(\beta,L/2)$.
Let $\lambda:\bR\to\bR$ be a symmetric function satisfying

\vskip0.2cm

\centerline{$\lambda\in\bH(\beta,1/2),\quad\int_{\bR}\lambda(y)\rd y=0,\quad
\lambda(0)>0,\quad \lambda(y)=0,\;\forall y\notin[-1,1].
$}

\vskip0.1cm

\noindent Let $\varpi>0$ be a constant the choice of which will be done later. Set  $\e=\big(\varpi L^{-2}\mn(\cQ_{\delta_n})/n\big)^{\frac{1}{2\beta+1}}$ and let

\vskip0.2cm

\centerline{$\rp(y)=\mathbf{n}(y)+L\e^\beta\lambda(y/\e),\quad y\in\bR.
$}

\vskip0.2cm

\noindent Obviously, $\rp\in\bH(\beta,L)$, $\int_{\bR}\rp(y)\rd y=1$ and for all $n$ large enough $\rp>0$. Hence $\rp$ is a probability density. Define

\vskip0.2cm

\centerline{$g(x)=\rp(x_1)\rp(x_2),\quad \mathbf{N}(x)=\mathbf{n}(x_1)\mathbf{n}(x_2),\; x=(x_1,x_2)\in\bR^2.
$}

\vskip0.2cm

\noindent We can assert that $g,\mathbf{N}\in\cG(\beta,L)$.
Let $\cQ_{\delta_n}=\{Q_0,\ldots Q_{\mm_n}\}$, where we have denoted $\mm_n=\exp\{\mn(\cQ_{\delta_n})\}-1=\text{card}(\cQ_{\delta_n})-1$. Set finally
\vskip0.2cm

\centerline{$f_0(\bullet)=\mathbf{N}(\bullet),\quad f_j(\bullet)=g(Q_j^T\bullet),\;\; j=1,\ldots,\mm_n.
$}

\vskip0.2cm

\noindent We assert that $\{f_j, j=0,\ldots,\mm_n\}\subset\cF\big(\beta,L,\cQ_{\delta_n}\big)$. Here we have used that $\mathbf{N}(\bullet)\equiv\mathbf{N}(Q_0^T\bullet)$. Additionally, for any $j=1,\ldots,\mm_n$
\begin{equation}
\label{eq1:lb}
\big|f_j(0)-f_0(0)|\big|=\big(\mathbf{n}(0)+L\lambda(0)\e^{\beta}\big)^{2}-\mathbf{n}^2(0)\geq c\big(\varpi L^{\frac{1}{\beta}}\mn(\cQ_{\delta_n})/n\big)^{\frac{\beta}{2\beta+1}},
\end{equation}
for all $n$ large enough. Here $c>0$ is a numerical constant independent on $n$ and $L$. Introduce

\vskip0.0cm

\centerline{$\displaystyle{Z_n=\frac{1}{\mm_n}\sum_{j=1}^{\mm_n}\prod_{i=1}^n\frac{f_j(X_i)}{f_0(X_i)}}.
$}

\vskip0.0cm

\noindent In view of (\ref{eq1:lb}) and in accordance with Corollary 2 of Proposition 5 in \cite{klp2} the assertion of the theorem will follow
with
$ \mathbf{c}_2=2^{-p}\Big(1-\sqrt{\frac{\Upsilon+1}{\Upsilon+5}}\Big)$
if we prove that
\begin{equation}
\label{eq2:lb}
\Upsilon:=\limsup_{n\to\infty}\bE_{f_0}\big\{Z_n^2\big\}<\infty.
\end{equation}
\quad $\mathbf{2^0}.\;$ We have
\begin{gather}
\bE_{f_0}\big\{Z_n^2\big\}=\frac{1}{\mm^2_n}\sum_{j=1}^{\mm_n}\bigg(\int_{\bR^2}\frac{f^2_j(x)}{f_0}\rd x\bigg)^{n}+
\frac{1}{\mm^2_n}\sum_{\stackrel{k,j=1}{k\neq j}}^{\mm_n}\bigg(\int_{\bR^2}\frac{f_j(x)f_k(x)}{f_0}\rd x\bigg)^{n}
\nonumber\\
\label{eq3:lb}
\leq \mm_n^{-1}\sup_{j=1,\ldots,\mm_n}\bigg(\int_{\bR^2}\frac{f^2_j(x)}{f_0}\rd x\bigg)^{n}+\sup_{\stackrel{k,j=1,\ldots,\mm_n}{k\neq j}}\bigg(\int_{\bR^2}\frac{f_j(x)f_k(x)}{f_0}\rd x\bigg)^{n}.
\end{gather}
\quad $\mathbf{2^0a}.\;$ Denote by $\mathbf{M}_{\e}(\cdot)=\mathbf{M}_1(\cdot)+\mathbf{M}_2(\cdot)+L\e^\beta\Lambda(\cdot)$, where we put

\vskip0.2cm

\centerline{$ \mathbf{M}_1(x)=\mathbf{n}(x_1)\lambda(x_2/\e),\quad \mathbf{M}_2(x)=\lambda(x_1/\e)\mathbf{n}(x_2),\quad\Lambda=\lambda(x_1/\e)\lambda(x_2/\e).
$}

\vskip0.2cm

\noindent Since $\int\lambda=0$ we have
\begin{equation}
\label{eq33:lb}
\int_{\bR^2}\mathbf{M}_\e\big(Q_j^Tx\big)\rd x=0,\quad \forall j=1,\ldots,\mm_n.
\end{equation}
Note also that for any $k,j=1,\ldots,\mm_n$
\begin{gather*}
f_j(x)f_k(x)=[\mathbf{N}(x)+L\e^{\beta}\mathbf{M}_\e\big(Q_j^Tx\big)][\mathbf{N}(x)+L\e^{\beta}\mathbf{M}_\e\big(Q_k^Tx\big)]
\\
=\mathbf{N}^2(x)+L\e^\beta\mathbf{N}(x)\big[\mathbf{M}_\e\big(Q_j^Tx\big)+\mathbf{M}_\e\big(Q_k^Tx\big)\big]+L^2\e^{2\beta}
\mathbf{M}_\e\big(Q_j^Tx\big)\mathbf{M}_\e\big(Q_k^Tx\big).
\end{gather*}
Thus, in view of (\ref{eq33:lb}) we have for any $j,k=1,\ldots,\mm_n$
\begin{gather}
\label{eq:0000}
a_{j,k}:=\int_{\bR^2}\frac{f_j(x)f_k(x)}{f_0(x)}\rd x=1+L^2\e^{2\beta}\int_{\bR^2}\frac{\mathbf{M}_\e\big(Q_j^Tx\big)\mathbf{M}_\e\big(Q_k^Tx\big)}{\mathbf{N}(x)}\rd x.
\end{gather}
\quad $\mathbf{2^0b}.\;$ It yields first,
\begin{eqnarray*}
a_{j,j}&=&1+L^2\e^{2\beta}\int_{\bR^2}\frac{\mathbf{M}_\e^2(x)}{\mathbf{N}(x)}\rd x
\\
&\leq&
1+3L^2\e^{2\beta}\bigg[\int_{\bR^2}\frac{\mathbf{M}_1^2(x)}{\mathbf{N}(x)}\rd x+\int_{\bR^2}\frac{\mathbf{M}_2^2(x)}{\mathbf{N}(x)}\rd x+L^{2}\e^{2\beta}\int_{\bR^2}\frac{\Lambda^2(x)}{\mathbf{N}(x)}\rd x\bigg]
\\
&=&1+3L^2\e^{2\beta}\bigg[2\int_{\bR}\frac{\lambda^2(y/\e)}{\mathbf{n}(y)}\rd y+
L^{2}\e^{2\beta}\bigg(\int_{\bR}\frac{\lambda^2(y/\e)}{\mathbf{n}(y)}\rd y\bigg)^2\bigg].
\end{eqnarray*}
From now on we will assume that $n$ is sufficiently large to guarantee that $\mathbf{n}(y)\geq 2^{-1}\mathbf{n}(0)$ for all $y\in [-\e,\e]$. Then, taking into account that $\lambda(y/\e)=0$ for any $y\notin [-\e,\e]$  we obtain for all $n$ large enough

\vskip0.15cm

\centerline{$ a_{j,j}
\leq 1+C_1L^2\e^{2\beta+1}.
$}

\vskip0.1cm

\noindent where $C_1$ is independent on $n$ and $L$. Hence, choosing $\varpi=C_1^{-1}$  we get
\begin{gather}
\label{eq4:lb}
a^n_{j,j}\leq e^{nC_1L^2\e^{2\beta+1}}=e^{C_1\varpi\mn(\cQ_{\delta_n})}=\mm_n,\quad \forall j=1,\ldots,\mm_n.
\end{gather}
\quad $\mathbf{2^0c}.\;$ For any $j,k=1,\ldots,\mm_n$, $j\neq k$ introduce $\cP_{j,k}=Q_j^TQ_k$. We have
\begin{eqnarray*}
b_{j,k}&:=&\int_{\bR^2}\frac{\mathbf{M}_\e\big(Q_j^Tx\big)\mathbf{M}_\e\big(Q_k^Tx\big)}{\mathbf{N}(x)}\rd x=\int_{\bR^2}\frac{\mathbf{M}_\e(x)\mathbf{M}_\e\big(\cP_{j,k}x\big)}{\mathbf{N}(x)}\rd x
\\
&=&
\int_{\bR^2}\frac{\lambda(x_2/\e)\mathbf{M}_\e\big(\cP_{j,k}x\big)}{\mathbf{n}(x_2)}\rd x+\int_{\bR^2}\frac{\lambda(x_1/\e)\mathbf{M}_\e\big(\cP_{j,k}x\big)}{\mathbf{n}(x_1)}\rd x
\\
&&\;+\int_{\bR^2}\frac{\Lambda(x)\mathbf{M}_\e\big(\cP_{j,k}x\big)}{\mathbf{N}(x)}\rd x.
\end{eqnarray*}
Taking into account that $\Lambda(x)=0$ for any $x\notin [-\e,\e]^2$, $\mathbf{N}(x)\geq4^{-1}\mathbf{n}^2(0)$ on $[-\e,\e]^2$ and $M_\e$
is uniformly bounded, we obtain that  for all $n$ large enough and some $C_2$ independent on $n$ and $L$
\begin{eqnarray}
\label{eq5:lb}
&&\int_{\bR^2}\frac{|\Lambda(x)|\big|\mathbf{M}_\e\big(\cP_{j,k}x\big)\big|}{\mathbf{N}(x)}\rd x\leq C_2L\e^{2}.
\end{eqnarray}
Also, we have for sufficiently large $n$
\begin{eqnarray*}
&&\int_{\bR^2}\frac{\lambda(x_2/\e)\mathbf{M}_\e\big(\cP_{j,k}x\big)}{\mathbf{n}(x_2)}\rd x\leq 2\mathbf{n}^{-1}(0)\bigg[
\int_{\bR^2}\lambda(x_2/\e)\mathbf{M}_1\big(\cP_{j,k}x\big)\rd x
\\
&&+\int_{\bR^2}\lambda(x_2/\e)\mathbf{M}_2\big(\cP_{j,k}x\big)\rd x+L\e^{\beta}\int_{\bR^2}\lambda(x_2/\e)\Lambda\big(\cP_{j,k}x\big)\rd x\bigg].
\end{eqnarray*}
Putting for brevity $\bp_1=p_1(Q_j,Q_k)$ and $\bp_2=p_2(Q_j,Q_k)$ and making the change of variables: $\e z_1=\bp_1x_1+\bp_2x_2, \e z_2= x_2$ (first and third integrals), $\e z_1=\bp_2x_1-\bp_1z_2, \e z_2= x_2$ (second integral) we obtain  since $\mathbf{n}$ and $\lambda$ are uniformly bounded for all $n$ large enough
\begin{eqnarray*}
&&\bigg|\int_{\bR^2}\frac{\lambda(x_2/\e)\mathbf{M}_\e\big(\cP_{j,k}x\big)}{\mathbf{n}(x_2)}\rd x\bigg|\leq \big[|\bp_1|\wedge|\bp_2|\big]^{-1}
\big[C_3\e^2+C_4L\e^{2+\beta}\big],
\end{eqnarray*}
where $C_3$ and $C_4$ are the constants independent of $n$ and $L$.

Since $|\bp_1|\wedge|\bp_2|\geq \delta_n$ in view of the definition of $\cQ_{\delta_n}$ we obtain
\begin{eqnarray}
\label{eq6:lb}
&&\bigg|\int_{\bR^2}\frac{\lambda(x_2/\e)\mathbf{M}_\e\big(\cP_{j,k}x\big)}{\mathbf{n}(x_2)}\rd x\bigg|\leq
C_5\e^2\delta_n^{-1}.
\end{eqnarray}
By the same computation we get
\begin{eqnarray}
\label{eq7:lb}
&&\bigg|\int_{\bR^2}\frac{\lambda(x_1/\e)\mathbf{M}_\e\big(\cP_{j,k}x\big)}{\mathbf{n}(x_1)}\rd x\bigg|\leq
C_5\e^2\delta_n^{-1}.
\end{eqnarray}
Collecting the bounds obtained in (\ref{eq5:lb}),  (\ref{eq6:lb}) and  (\ref{eq7:lb}) we obtain

\vskip0.15cm

\centerline{$
 \big|b_{j,k}\big|\leq  C_6\e^2\delta_n^{-1}\leq C_7\e\mn^{-1}(\cQ_{\delta_n}),\quad\forall j,k=1,\;\ldots,\mm_n,\; j\neq k.
$}

\vskip0.15cm

\noindent
because $\delta_n\geq \big(\ln(n)\big)^{\frac{2\beta+2}{2\beta+1}}\big(L^{-2}/n\big)^{\frac{1}{2\beta+1}}$ in view of the assumption of the theorem and $\mn(\cQ_{\delta_n})\leq C_8\ln(n)$ in view of (\ref{eq:bound-n_delta}). It yields together with (\ref{eq:0000})

\vskip0.15cm

\centerline{$
\big|a_{j,k}\big|^n\leq e^{C_9}
$}

\vskip0.15cm

\noindent for any $j,k=1,\;\ldots,\mm_n, j\neq k$ and all $n$ large enough. This in its turn, together with (\ref{eq4:lb}) and (\ref{eq3:lb}) allows us to assert that (\ref{eq2:lb}) holds with
$\Upsilon\leq 1+e^{C_9}$. The proof of the theorem is completed.
\epr

\subsection{Proofs of Theorems \ref{th2} and \ref{th3}}

The proofs of Theorems \ref{th2} and \ref{th3} are essentially based on several auxiliary results. We starts with presenting such of them which
will be used in the proofs of the both theorems simultaneously. Their proofs as well as the proofs of all auxiliary results are postponed to Appendix section.
Set  for any $f\in\cF(\beta,L,\mQ)$, $D\in\mQ$ and $x\in\bR^2$

\vskip0.25cm

\centerline{$
\tau_f(D)=\left\{
\begin{array}{cc}
\int_{\bR^2}\Bg_f\big(\bp_1\Gamma u\big) \Bg_f\big(\bp_1^{-1}D\Omega x+\bp_2\Omega\Gamma u\big)\rd u,\quad &D\neq\BQ_f;
\\*[2mm]
f(x),\quad &D=\BQ_f.
\end{array}
\right.
$}

\vskip0.2cm

\noindent
where we have put
$\bp_1=p_1(D,\BQ_f), \bp_2=p_2(D,\BQ_f)$. Set also

\vskip0.1cm

\centerline{$
\displaystyle{\widetilde{f}_{h,d}(x)=n^{-1}\sum_{k=1}^n\cK_{h}\big(d^T(X_k-x)\big),\;
\widetilde{f}_{h,d_\perp}(x)=n^{-1}\sum_{k=1}^n\cK_{h}\big(d^T_\perp(X_k-x)\big)}.
$}

\vskip-0.3cm

\begin{lemma}
\label{lem1}
For any  $D\in\mQ$, $\beta>0$, $L>0$,  $x\in\bR^2$, $\cK\in\bK_\beta$ and  $h>0$
\vskip0.2cm

\centerline{$
\displaystyle{\sup_{f\in\cF(\beta,L,\mQ)}\Big|\bE_{f}\big[\widetilde{f}_{h,d}(x)\big]\bE_{f}\big[\widetilde{f}_{h,d_\perp}(x)\big]-\tau_f\big(D\big)\Big|
\leq2C(\cK,\beta,1)L^2 h^\beta}.
$}

\end{lemma}


\begin{lemma}
\label{lem:new1}
For any  $\beta>0$, $L>0$, $x\in\bR^2$ and $\cK\in\bK_\beta$
\begin{eqnarray*}
&&\sup_{D\in\mQ}\sup_{f\in\cF(\beta,L,\mQ)}\Big|\bE_{f}\big[\widetilde{f}_{h,(D,\BQ_f)}(x)\big]-
\bE_{f}\big[\widetilde{f}_{\eta,d}(x)\big]\bE_{f}\big[\widetilde{f}_{\eta,d_\perp}(x)\big]\Big|
\\
&&\qquad\quad\quad\;\;\leq 2C\big(\cK,\beta,\sqrt{2}\big)L^2 \big(h^\beta+\eta^\beta\big),\quad\forall h,\eta>0.
\end{eqnarray*}

\end{lemma}


\begin{lemma}
\label{lem:new2}
For any   $D,Q\in\mQ,$ and any  $f\in\cF(\beta,L,\mQ)$

\vskip0.1cm

\centerline{$
\displaystyle{\bE_{f}\Big[\widetilde{f}_{h,(D,Q)}(x)\Big]=\bE_{f}\Big[\widetilde{f}_{h,(Q,D)}(x)\Big]}.
$}

\end{lemma}

\vskip-0.2cm

\noindent This feature of the auxiliary estimator was called in \cite{GL12} the \textsf{commutativity property}.

Let $\mI_n$ be the set of all pairwise disjoint subsets of $\{1,\ldots,n\}$. For any $\cI\in\mI_n$ its cardinality is denoted by $|\cI|$ and
$\widetilde{f}^{(\cI)}_{h,Q}$ will be used for the estimator built from $(X_i, i\in\cI)$.

\begin{proposition}
\label{prop2}
Let $p\geq 1$,  $\beta>0$, $L>0$,   and  $\cK\in\bK_\beta$  be fixed and set $\blh=\big(\mu/|\cI|\big)^{\frac{1}{2\beta+1}}, \mu>0,\cI\in\mI_n$.
There exist $\mathbf{c}_3$ independent of $L$ such that

\vskip0.1cm

\centerline{$
 \displaystyle{\sup_{n\geq 1}\sup_{\cI\in\mI_n}\sup_{\mu\in[1,|\cI|]}(\mu/|\cI|)^{-\frac{\beta}{2\beta+1}}\sup_{\cF(\beta,L,\mQ)} \cR^{(p)}_n\big[\widetilde{f}^{(\cI)}_{\blh,\BQ_f}, f\big]\leq \mathbf{c}_3L\big(L+L^{\frac{1}{p\vee 2}})}.
 $}

 \end{proposition}

\subsubsection{Proof of Theorem \ref{th2}}
\label{th2:proof}
Let us formulate some auxiliary results the proofs of which are postponed to Appendix section.

\vskip0.1cm

Put $\mathbf{A}=12\sqrt{10p\alpha}\left(1+\sqrt{5p}\right)\left[1\vee\|\cK\|_\infty\right]+4C(\cK,\mb,\sqrt{2})$, where

\vskip0.1cm

\noindent $\displaystyle{\alpha=1\vee\sup_{n\geq 3}\left\{[1\vee\mn(\cQ_{\delta_n})]/\ln(n)\right\}}$ is finite in view of (\ref{eq:bound-n_delta}) and (\ref{eq:delta}).

\vskip0.1cm

\noindent Set for any  $n\geq 3$,  $\delta_n>0$,  $\cQ_{\delta_n}\in\bQ_{\delta_n}$  and $f\in\cF(\beta,L,\cQ_{\delta_n})$
\begin{gather}
\zeta_n(f,x)=\sup_{\stackrel{h\in\mH,}{D,Q\in\cQ_{\delta_n}}}\Big[\big|\widetilde{f}_{h,(D,Q)}(x)-\varkappa_h(D,Q,x)
\big|-\boldsymbol{\ma}\widehat{\cU}_n\sqrt{\ln(n)/nh}\Big]_+.
\nonumber\\
\label{eq0:proof}
\hskip-0cm \varkappa_h(D,Q,x)=\left\{
\begin{array}{ll}
\bE_{f}\big[\widetilde{f}_{h,(D,Q)}(x)\big],\quad&D\neq Q;
\\*[1mm]
\bE_{f}\big[\widetilde{f}_{h,q}(x)\big]\bE_{f}\big[\widetilde{f}_{h,q_\perp}(x)\big],\quad&D=Q,
\end{array}
\right.
\end{gather}
where $\boldsymbol{\ma}=2^{-1}\mathbf{A}-2C\big(\cK,\mb,\sqrt{2}\big)$.

\begin{proposition}
\label{prop1}
For any uniformly bounded kernel $\cK$, an arbitrary sequence $\delta_n$ satisfying (\ref{eq:delta}) and
 any $\cQ_{\delta_n}\in\bQ_{\delta_n}$ one has
\begin{eqnarray*}
&&\limsup_{n\to\infty}n^{3p}\sup_{\beta\in\{\beta_1,\beta_2\}}\sup_{f\in\cF(\beta,L,\cQ_{\delta_n})}\bP_f\big(\zeta_n(f,x)\neq 0\big)=0;
\\
&&\limsup_{n\to\infty}n^{p}\sup_{\beta\in\{\beta_1,\beta_2\}}\sup_{f\in\cF(\beta,L,\cQ_{\delta_n})}\bE_{f}\big[\zeta^p_n(f,x)\big]=0;
\\*[1mm]
&&\limsup_{n\to\infty}\sup_{\beta\in\{\beta_1,\beta_2\}}\sup_{f\in\cF(\beta,L,\cQ_{\delta_n})}\bE_{f}\big(\widehat{\cU}^p_n\big)\;\leq\; \mathbf{C}_p(\cK)L^{2p},
\end{eqnarray*}
where $\mathbf{C}_p(\cK)>0$ is given in the proof of the proposition.

\end{proposition}

\noindent {\it Proof of the theorem.\;} We divide the proof into several steps.

$\mathbf{1^0}.\;$  For any $\beta\in (0,\mb]$ and $L>0$ set $\mh=\big(L^{-4}\ln(n)/n\big)^{\frac{1}{2\beta+1}}$ and let $\cA=\big\{R_{n}(\BQ_f,\mh)\neq 0\big\}$.
Our first goal is to prove the following result.
\begin{eqnarray}
\label{eq1:proof}
\lim_{n\to\infty}\sup_{\beta\in\{\beta_1,\beta_2\}}\big(L^{\frac{2}{\beta}}\ln(n)/n\big)^{-\frac{p\beta}{2\beta+1}}\sup_{\cF(\beta,L,\cQ_{\delta_n})} \bE_f\big|\widehat{f}- f(x)\big|^p\mathrm{1}_{\cA}=0.
\end{eqnarray}
Note that for any  $n\geq 1$ since $\widehat{h}\in\cH$

\vskip0.15cm

\centerline{$
 \displaystyle{\big|\widetilde{f}_{\widehat{h},\widehat{Q}}(x)\big|\leq \|\cK\|^2_{\infty}\widehat{h}^{-2}\leq \|\cK\|^2_{\infty}n^2}.
 $}

\vskip0.15cm

\noindent Hence (\ref{eq1:proof}) will be proved if we show that
\begin{eqnarray}
\label{eq001:proof}
&&\quad\limsup_{n\to\infty}\sup_{\beta\in\{\beta_1,\beta_2\}}n^{3p}\sup_{\cF(\beta,L,\cQ_{\delta_n})} \bP_f(\cA)=0.
\end{eqnarray}
For any $\eta,\eta^\prime\in\mH$, $\eta^\prime\leq\eta\leq \mh$ we have in view of the definition of $\varkappa_\cdot(\cdot,\cdot,x)$
\begin{gather*}
\label{eq01:proof}
\sup_{D\in\cQ_{\delta_n}}\big|\widetilde{f}_{\eta,(D,\BQ_f)}(x)-\widetilde{f}_{\eta^\prime,D}(x)\big|\leq \sup_{D\in\cQ_{\delta_n}}\big|\widetilde{f}_{\eta,(D,\BQ_f)}(x)-\varkappa_\eta(D,\BQ_f,x)\big|
\\
\quad+\sup_{D\in\cQ_{\delta_n}}\big|\widetilde{f}_{\eta^\prime,D}(x)-\varkappa_{\eta^\prime}(D,D,x)\big|+
\sup_{D\in\cQ_{\delta_n}}\big|\varkappa_\eta(D,\BQ_f,x)-\varkappa_{\eta^\prime}(D,D,x)\big|
\\
\hskip-5.9cm\leq 2{\boldsymbol\ma}\widehat{\cU}_n
\big(\ln(n)/n\eta^\prime\big)^{1/2}+2\zeta_n(f,x)
\\
\hskip-2.4cm +\sup_{D\in\cQ_{\delta_n}}\big|\bE_{f}\big[\widetilde{f}_{\eta,(D,\BQ_f)}(x)\big]-
\bE_{f}\big[\widetilde{f}_{\eta^\prime,d}(x)\big]\bE_{f}\big[\widetilde{f}_{\eta^\prime,d_\perp}(x)\big]\big|
\end{gather*}
Taking into account that $(\ln(n)/n\mh)^{1/2}=L^2\mh^\beta$, $\widehat{\cU}_n\geq 1$ and applying  Lemma \ref{lem:new1} we get in view of the definition of $\mathbf{A}$ for any $\eta^\prime\leq\eta\leq \mh$
\begin{gather*}
\sup_{D\in\cQ_{\delta_n}}\big|\widetilde{f}_{\eta,(D,\BQ_f)}(x)-\widetilde{f}_{\eta^\prime,D}(x)\big|\leq 2{\boldsymbol\ma}\widehat{\cU}_n
\big(\ln(n)/n\eta^\prime\big)^{1/2}+2\zeta_n(f,x)
\\
+2L^2C\big(K,\mb,\sqrt{2}\big)\big]\big(\eta^{\beta}+(\eta^\prime)^\beta\big)\leq \widehat{\cU}_n\big[2{\boldsymbol\ma}+4C\big(\cK,\mb,\sqrt{2}\big)\big]\big(\ln(n)/n\eta^\prime\big)^{1/2}
\nonumber\\*[1mm]
+2\zeta_n(f,x)=\mathbf{A}\widehat{\cU}_n\big(\ln(n)/n\eta^\prime\big)^{1/2}+2\zeta_n(f,x).
\end{gather*}
\noindent Thus we have

\vskip0.1cm

\centerline{$
\displaystyle{\sup_{\stackrel{\eta,\eta^\prime\in\mH:}{\eta^\prime\leq\eta\leq \mh}}\sup_{D\in\cQ_\delta }\big[\big|\widetilde{f}_{\eta,(D,\BQ_f)}(x)-\widetilde{f}_{\eta^\prime,D}(x)\big|-\mathbf{A}\widehat{\cU}_n
\big(\ln(n)/n\eta^\prime\big)^{1/2}
\;\Big]_+\leq 2\zeta_n}.
$}

\vskip0.1cm

\noindent
The latter means that

\vskip0.1cm

\centerline{$
\bP_f(\cA)\leq \bP_f\big(\zeta_n(f,x)\neq 0\big)
$}

\vskip0.1cm

\noindent and (\ref{eq001:proof}) follows from the first assertion of Proposition \ref{prop1}.

 $\mathbf{2^0}.\;$ Denote  $\bar{\cA}$ the event complimentary to $\cA$.  Note  that if $\bar{\cA}$ is realized
\begin{eqnarray}
\label{eq2:proof}
R_{n}\big(\widehat{Q},\widehat{h}\big)&\leq& R_{n}\big(\widehat{Q},\widehat{h}\big)+\mathbf{A}\widehat{\cU}_n\sqrt{\ln(n)/n\widehat{h}}
\nonumber\\
&\leq&  R_{n}\big(\BQ_f,\mh\big)+\mathbf{A}\widehat{\cU}_n\sqrt{\ln(n)/n\mh}=\mathbf{A}\widehat{\cU}_n\sqrt{\ln(n)/n\mh}.
\end{eqnarray}
To get the second inequality we have used
the definition of $(\widehat{h},\widehat{Q})$.

 Now let us prove the following inclusion.
\begin{eqnarray}
\label{eq3:proof}
&& \bar{\cA}\subseteq \{\widehat{h}\geq \mh\},
\end{eqnarray}
Indeed, if $\bar{\cA}$ is realized then
\begin{eqnarray*}
\mathbf{A}\widehat{\cU}_n\sqrt{\ln(n)/n\widehat{h}}&\leq& R_{n}\big(\widehat{Q},\widehat{h}\big)+\mathbf{A}\widehat{\cU}_n\sqrt{\ln(n)/n\widehat{h}}
\\
&\leq&  R_{n}\big(\BQ_f,\mh\big)+\mathbf{A}\widehat{\cU}_n\sqrt{\ln(n)/n\mh}=\mathbf{A}\widehat{\cU}_n\sqrt{\ln(n)/n\mh}
\end{eqnarray*}
and (\ref{eq3:proof}) follows.

 $\mathbf{3^0}.\;$ If $\bar{\cA}$ is realized and, therefore $\widehat{h}\geq \mh$ in view of (\ref{eq3:proof}),  we have

\vskip0.15cm

\centerline{$\Big|\widetilde{f}_{\widehat{h},\widehat{Q}}(x)-\widetilde{f}_{\mh,\widehat{Q}}(x)\Big|\leq R_{n}\big(\widehat{Q},\widehat{h}\big)+\mathbf{A}\widehat{\cU}_n\sqrt{\ln(n)/n\mh}.
$}

\vskip0.15cm

\noindent
This yields together with (\ref{eq2:proof})
\begin{eqnarray}
\label{eq4:proof}
&& \Big|\widetilde{f}_{\widehat{h},\widehat{Q}}(x)-\widetilde{f}_{\mh,\widehat{Q}}(x)\Big|\mathrm{1}_{\bar{\cA}}\leq 2\mathbf{A}\widehat{\cU}_n\sqrt{\ln(n)/n\mh}.
\end{eqnarray}
Putting $\cB=\{\widehat{Q}=\BQ_f\}$ we deduce from (\ref{eq4:proof})
\begin{eqnarray}
\label{eq5:proof}
&& \Big|\widetilde{f}_{\widehat{h},\widehat{Q}}(x)-f(x)\Big|\mathrm{1}_{\bar{\cA}\cap\cB}\leq 2\mathbf{A}\widehat{\cU}_n\sqrt{\ln(n)/n\mh}+
\big|\widetilde{f}_{\mh,\BQ_f}(x)-f(x)\big|.
\end{eqnarray}
Also we obtain using (\ref{eq2:proof})
\begin{eqnarray}
\label{eq6:proof}
\hskip-2cm\big|\widetilde{f}_{\mh,(\BQ_f,\widehat{Q})}(x)-\widetilde{f}_{\mh,\BQ_f}(x)\big|\mathrm{1}_{\bar{\cA}\cap\bar{\cB}}&\leq&
R_{n}\big(\widehat{Q},\widehat{h}\big)+\mathbf{A}\widehat{\cU}_n\sqrt{\ln(n)/n\mh}
\\
&\leq& 2\mathbf{A}\widehat{\cU}_n\sqrt{\ln(n)/n\mh};
\nonumber
\\*[1mm]
\label{eq7:proof}
\big|\widetilde{f}_{\mh,(\widehat{Q},\BQ_f)}(x)-\widetilde{f}_{\mh,\widehat{Q}}(x)\big|\mathrm{1}_{\bar{\cA}\cap\bar{\cB}}&\leq&
R_{n}\big(\BQ_f,\mh\big)+\mathbf{A}\widehat{\cU}_n\sqrt{\ln(n)/n\mh}
\\
&=& \mathbf{A}\widehat{\cU}_n\sqrt{\ln(n)/n\mh}.
\nonumber
\end{eqnarray}
 We have in view of Lemma \ref{lem:new2}
\begin{gather}
\label{eq8:proof}
\big|\widetilde{f}_{\mh,(\BQ_f,\widehat{Q})}(x)-\widetilde{f}_{\mh,(\widehat{Q},\BQ_f)}(x)\big|\leq
\sup_{D,Q\in\cQ_{\delta_n}}\big|\widetilde{f}_{\mh,(Q,D)}(x)-\widetilde{f}_{\mh,(D,Q)}(x)\big|
\\
\leq\sup_{D,Q\in\cQ_{\delta_n}}
\big|\bE_f\big[\widetilde{f}_{\mh,(Q,D)}(x)\big]-\bE_f\big[\widetilde{f}_{\mh,(D,Q)}(x)\big]\big|
\nonumber
\\*[1mm]
\quad\;+
2\mathbf{A}\widehat{\cU}_n\sqrt{\ln(n)/n\mh}+2\zeta_n(f,x)=2\mathbf{A}\widehat{\cU}_n\sqrt{\ln(n)/n\mh}+2\zeta_n(f,x).
\nonumber
\end{gather}
\quad$\mathbf{4^0}.\;$ We obtain from (\ref{eq4:proof}), (\ref{eq6:proof}), (\ref{eq7:proof}) and (\ref{eq8:proof})

\vskip0.15cm

\centerline{$\Big|\widetilde{f}_{\widehat{h},\widehat{Q}}(x)-f(x)\Big|\mathrm{1}_{\bar{\cA}\cap\bar{\cB}}\leq 7\mathbf{A}\widehat{\cU}_n\sqrt{\ln(n)/n\mh}+2\zeta_n(f,x)+
\big|\widetilde{f}_{\mh,\BQ_f}(x)-f(x)\big|.
$}

\vskip0.15cm

\noindent It yields together with (\ref{eq5:proof})
\begin{equation}
\label{eq9:proof}
\Big|\widetilde{f}_{\widehat{h},\widehat{Q}}(x)-f(x)\Big|\mathrm{1}_{\bar{\cA}}\leq 7\mathbf{A}\widehat{\cU}_n\sqrt{\ln(n)/n\mh}+2\zeta_n(f,x)+
\big|\widetilde{f}_{\mh,\BQ_f}(x)-f(x)\big|.
\end{equation}
 Since $\ln(n)/n\mh=(L^{\frac{2}{\beta}}\ln(n)/n)^{\frac{2\beta}{2\beta+1}}$, we deduce from the second and third assertions of  Proposition \ref{prop1}, (\ref{eq1:proof}) and
(\ref{eq9:proof})
\begin{gather}
\limsup_{n\to\infty}\sup_{\beta\in\{\beta_1,\beta_2\}}\big(L^{\frac{1}{\beta}}\ln(n)/n\big)^{-\frac{\beta}{2\beta+1}}\sup_{\cF(\beta,L,\cQ_{\delta_n})}
\cR^{(p)}_n\big[\widehat{f}, f\big]\leq C_p\Big\{L^2
\nonumber\\
\label{eq10:proof}
+
\limsup_{n\to\infty}\sup_{\beta\in\{\beta_1,\beta_2\}}\big(L^{\frac{1}{\beta}}\ln(n)/n\big)^{-\frac{\beta}{2\beta+1}}\sup_{\cF(\beta,L,\cQ_{\delta_n})}
\cR^{(p)}_n\big[\widetilde{f}_{\mh,\BQ_f}(x), f\big]\Big\},
\end{gather}
where $C_p$ depends on $p$  and $\cK$ only.

The assertion of the theorem follows now from  (\ref{eq10:proof}) and Proposition \ref{prop2} where one should choose  $\mu=\min_{\beta\in\{\beta_1,\beta_2\}}L^{2/\beta}\ln(n)$ and $\cI=\{1,\ldots,n\}$.
\epr

\subsubsection{Proof of Theorem \ref{th3}}
\label{th3:proof} The proof of the theorem is similar to those of Theorem \ref{th2} and essentially based on the  result formulated in  Proposition \ref{prop4} below.

\vskip0.1cm

Put $\mathbf{B}=527730p^2\sqrt{6}\left(\|\cK\|_1^2\vee\|\cK\|_2^2\vee\|\cK\|_\infty^2\right)[9+4\alpha]^{\frac{3\beta+3}{2\beta+1}}\left[C(\beta)\right]^{\frac{3}{2}}L^{\frac{4\beta+8}{2\beta+1}}+8C(\cK,\mb,\sqrt{2})L^2$, where $\displaystyle{C(\beta):=1\vee\sup_{n\geq 3}\Big\{\big[\ln^2(n)/n\big]^{\frac{2\beta}{2\beta+1}}\big[\ln(n)\big]^{\frac{2}{2\beta+1}}\Big\}}$.

\vskip0.1cm

\noindent Set for any  $n\geq 3$,  $\delta_n>0$,  $\cQ_{\delta_n}\in\bQ_{\delta_n}$
and $f\in\cF(\beta,L,\cQ_{\delta_n})$

\vskip0.1cm

\centerline{$
\displaystyle{\chi_\mfi(f,x)=\sup_{D,Q\in\cQ_{\delta_n}}\Big[\big|\widetilde{f}^{(\mathbf{i})}_{h_{\mathbf{i}},(D,Q)}(x)-\varkappa_{h_{\mathbf{i}}}(D,Q,x)
\big|-\mathbf{C} L^2h_{\mathbf{i}}^{\beta}\Big]_+},
$}

\vskip0.1cm

\noindent where $\varkappa_{h}(\cdot,\cdot,x), h>0,$ is defined in (\ref{eq0:proof}) and $\mathbf{C}=2^{-1}\mathbf{B}-4C\big(\cK,\beta,\sqrt{2}\big)L^2$.
\begin{proposition}
\label{prop4}
For any   $\beta>0, L>0$, $\cK\in\bK_\beta$, an arbitrary sequence $\delta_n$ satisfying (\ref{eq:delta}) with $\mb=\beta$,
 $\cQ_{\delta_n}\in\bQ_{\delta_n}$ and any $\mfi=1,\ldots,\mfi^*$ one has
\begin{gather*}
\sup_{n\geq 3} \sup_{\mfi=1,\ldots,\mfi^*} \bigg(\frac{\omega_{\mfi-1}\bn}{\mathbf{n}_{\mfi-1}\omega_{\mfi}}\bigg)^{\frac{p\beta}{2\beta+1}}
\sup_{f\in\cF(\beta,L,\cQ_{\delta_n})}\bP_{f}\big\{\chi_\mfi(f,x)\neq 0\big\}=:\mathbf{P}<1;
\\
\sup_{n\geq 3} \sup_{\mfi=1,\ldots,\mfi^*} (L^{\frac{2}{\beta}}\omega_{\mathbf{i}}/\bn)^{-\frac{p\beta}{2\beta+1}}\sup_{f\in\cF(\beta,L,\cQ_{\delta_n})}\bE_{f}\big\{\chi^p_\mfi(f,x)\big\}=:\mathbf{E}<\infty.
\end{gather*}

\end{proposition}

\noindent {\it Proof of the theorem.\;} Throughout the proof we will understand $\varkappa_{h_{\mathbf{i}}}(\cdot,\cdot,x)$ introduced in (\ref{eq0:proof}) as the mapping defined on $\cQ_{\delta_n}\times\cQ_{\delta_n}$ (its explicit expression via some integral operators can be easily obtained). It allows us to introduce below random variables $\varkappa_{h_{\mathbf{i}}}(\widehat{Q}^{(\mathbf{i})},\cdot,x), \mathbf{i}=1,\ldots,\mathbf{i}^*$.

$\mathbf{1^0}.\;$ Introduce the random event $\cZ^{(\mathbf{i})}=\big\{R^{(\mathbf{i})}_{n}(\widehat{Q}^{(\mathbf{i})})=0\big\}$. If $\cZ^{(\mathbf{i})}$ is realized, we assert,
using the definitions of $\widehat{Q}^{(\mathbf{i})}$ and $\breve{f}^{(\mathbf{i})}$ that
\begin{eqnarray}
\label{eq15:proof}
 R^{(\mathbf{i})}_{n}(\widehat{Q}^{(\mathbf{i})})=0\; &\Rightarrow&\;
\left\{
\begin{array}{ll}
\breve{f}^{(\mathbf{i})}(x)=\widetilde{f}^{(\mathbf{i})}_{h_\mfi, \widehat{Q}^{(\mathbf{i})}}(x);
\\*[2mm]
\Big|\widetilde{f}^{(\mathbf{i})}_{h_\mfi,(\BQ_f,\widehat{Q}^{(\mathbf{i})})}(x)-\widetilde{f}^{(\mathbf{i})}_{h_\mfi,\BQ_f}(x)\Big|\leq \mathbf{B} L^2h_{\mathbf{i}}^{\beta}.
\end{array}
\right.
\end{eqnarray}
Note that
 \begin{gather}
\big|\widetilde{f}^{(\mathbf{i})}_{h_\mfi, \widehat{Q}^{(\mathbf{i})}}(x)-f(x)\big|\leq
\big|\widetilde{f}^{(\mathbf{i})}_{h_\mfi, \widehat{Q}^{(\mathbf{i})}}(x)-\varkappa_{h_{\mathbf{i}}}(\widehat{Q}^{(\mathbf{i})},\widehat{Q}^{(\mathbf{i})},x)\big|
\nonumber\\
+
\big|\varkappa_{h_{\mathbf{i}}}(\widehat{Q}^{(\mathbf{i})},\widehat{Q}^{(\mathbf{i})},x)-f(x)\big|\leq \mathbf{C} L^2h_{\mathbf{i}}^{\beta}+\chi_\mfi(f,x)
\nonumber\\
+\sup_{Q\in\cQ_{\delta_n}}\big|\varkappa_{h_{\mathbf{i}}}(Q,Q,x)-\varkappa_{h_{\mathbf{i}}}(Q,\BQ_f,x)\big|+\big|\varkappa_{h_{\mathbf{i}}}
(\widehat{Q}^{(\mathbf{i})},\BQ_f,x)-f(x)\big|
\nonumber\\
\leq \mathbf{C} L^2h_{\mathbf{i}}^{\beta}+\chi_\mfi(f,x)+4C\big(\cK,\beta,\sqrt{2}\big)L^2h_{\mathbf{i}}^{\beta}+
\big|\varkappa_{h_{\mathbf{i}}}(\widehat{Q}^{(\mathbf{i})},\BQ_f,x)-f(x)\big|
\nonumber\\*[2mm]
\label{eq13:proof}
\leq 2^{-1}\mathbf{B} L^2h_{\mathbf{i}}^{\beta}+\chi_\mfi(f,x)+\big|\varkappa_{h_{\mathbf{i}}}(\widehat{Q}^{(\mathbf{i})},\BQ_f,x)-f(x)\big|.
\end{gather}
To get the penultimate inequality we have used Lemma \ref{lem:new1} while the last one follows from the definition of $\mathbf{C}$. Also in view of Lemma \ref{lem:new2} for all $D\in\cQ_{\delta_n}$
\begin{gather}
\label{eq13-new:proof}
\varkappa_{h_{\mathbf{i}}}(\widehat{Q}^{(\mathbf{i})},D,x)=\sum_{Q\in\cQ_{\delta_n}}\varkappa_{h_{\mathbf{i}}}(Q,D,x)
\mathrm{1}_{\widehat{Q}^{(\mathbf{i})}=Q}
\\
=\sum_{\stackrel{Q\in\cQ_{\delta_n}}{Q\neq D}}
\bE_f\big[\widetilde{f}^{(\mathbf{i})}_{h_\mfi,(Q,D)}(x)\big]\mathrm{1}_{\widehat{Q}^{(\mathbf{i})}=Q}
+\varkappa_{h_{\mathbf{i}}}(D,D,x)\mathrm{1}_{\widehat{Q}^{(\mathbf{i})}=D}
\nonumber\\
=\sum_{\stackrel{Q\in\cQ_{\delta_n}}{Q\neq D}}
\bE_f\big[\widetilde{f}^{(\mathbf{i})}_{h_\mfi,(D,Q)}(x)\big]\mathrm{1}_{\widehat{Q}^{(\mathbf{i})}=Q}
+\varkappa_{h_{\mathbf{i}}}(D,D,x)\mathrm{1}_{\widehat{Q}^{(\mathbf{i})}=D}=\varkappa_{h_{\mathbf{i}}}(D,\widehat{Q}^{(\mathbf{i})},x).
\nonumber
\end{gather}
Thus, if $\cZ^{(\mathbf{i})}$ is realized we have in view of (\ref{eq13-new:proof})
\begin{gather*}
\big|\varkappa_{h_{\mathbf{i}}}(\widehat{Q}^{(\mathbf{i})},\BQ_f,x)-f(x)\big|=\big|\varkappa_{h_{\mathbf{i}}}(\BQ_f,\widehat{Q}^{(\mathbf{i})},x)-f(x)\big|
\\
\leq \big|\varkappa_{h_{\mathbf{i}}}(\BQ_f,\widehat{Q}^{(\mathbf{i})},x)-\widetilde{f}^{(\mathbf{i})}_{h_\mfi,(\BQ_f,\widehat{Q}^{(\mathbf{i})})}(x)\big|+
\Big|\widetilde{f}^{(\mathbf{i})}_{h_\mfi,(\BQ_f,\widehat{Q}^{(\mathbf{i})})}(x)-\widetilde{f}^{(\mathbf{i})}_{h_\mfi,\BQ_f}(x)\Big|
\\
+
\Big|\widetilde{f}^{(\mathbf{i})}_{h_\mfi,\BQ_f}(x)-f(x)\Big|\leq \mathbf{C} L^2h_{\mathbf{i}}^{\beta}+\chi_\mfi(f,x)+\mathbf{B} L^2h_{\mathbf{i}}^{\beta}
+
\Big|\widetilde{f}^{(\mathbf{i})}_{h_\mfi,\BQ_f}(x)-f(x)\Big|.
\end{gather*}
It yields together with (\ref{eq13:proof})
\begin{gather}
\label{eq17:proof}
\big|\widetilde{f}^{(\mathbf{i})}_{h_\mfi, \widehat{Q}^{(\mathbf{i})}}(x)-f(x)\big|\mathrm{1}_{\cZ^{(\mathbf{i})}}
\leq 2\mathbf{B} L^2h_{\mathbf{i}}^{\beta}+2\chi_\mfi(f,x) +
\Big|\widetilde{f}^{(\mathbf{i})}_{h_\mfi,\BQ_f}(x)-f(x)\Big|.
\end{gather}
First, we deduce from Proposition \ref{prop2}

\vskip0.1cm

\centerline{$
\displaystyle{\sup_{n\geq 3}\sup_{\mfi=1,\ldots,\mathbf{i}^*}\big(L^{\frac{2}{\beta}}\omega_{\mathbf{i}}/\bn\big)^{-\frac{p\beta}{2\beta+1}}\sup_{\cF(\beta,L,\cQ_{\delta_n})}
\Big[\cR^{(p)}_n\big[\widetilde{f}^{(\mathbf{i})}_{h_\mfi,\BQ_f}, f\big]\Big]^p\Big\}=:C_1<\infty}.
$}

\vskip0.1cm

\noindent Next, taking into account that $L^2h_{\mathbf{i}}^{\beta}=(L^{\frac{2}{\beta}}\omega_{\mathbf{i}}/\bn)^{\frac{\beta}{2\beta+1}}$ and
 denoting

\vskip0.1cm

\centerline{$
\displaystyle{\mathbf{R}=\sup_{n\geq 3}\sup_{\mfi=1,\ldots,\mathbf{i}^*}\big(L^{\frac{2}{\beta}}\omega_{\mathbf{i}}/\bn\big)^{-\frac{p\beta}{2\beta+1}}\sup_{\cF(\beta,L,\cQ_{\delta_n})}
\bE_f\big|\breve{f}^{(\mathbf{i})}(x)-f(x)\big|^p\mathrm{1}_{\cZ^{(\mathbf{i})}}}
$}

\vskip0.1cm

\noindent
we deduce from (\ref{eq15:proof}), (\ref{eq17:proof}) and the second assertion of Proposition \ref{prop4}
\begin{eqnarray}
\label{eq18:proof}
&&\mathbf{R}
\leq 3^{p}\big\{(2\mathbf{B})^p+2^p\mathbf{E}+C_1 \big\}.
\end{eqnarray}
 \quad$\mathbf{2^0}.\;$ In view of the definition of $\breve{f}^{(\mathbf{i})}$ we have
 \begin{eqnarray}
\label{eq19:proof}
&&
\bE_f\big|\breve{f}^{(\mathbf{i})}(x)-f(x)\big|^p\mathrm{1}_{\bar{\cZ}^{(\mathbf{i})}}
=\bE_f\left\{\big|\breve{f}^{(\mathbf{i}-1)}(x)-f(x)\big|^p\mathrm{1}_{\bar{\cZ}^{(\mathbf{i})}}\right\}
\nonumber\\
&&=\bE_f\big|\breve{f}^{(\mathbf{i}-1)}(x)-f(x)\big|^p\;\bP_f\big(\bar{\cZ}^{(\mathbf{i})}\big),
\end{eqnarray}
since $X^{(\mathbf{i})}$ and $X^{(\mathbf{i}-1)}$ are the independent collections of random variables.
Note that in view of the definition of $\varkappa_{h_\mathbf{i}}(\cdot,\cdot,x)$
 \begin{gather*}
\sup_{D\in\cQ_{\delta_n}}\Big|\widetilde{f}^{(\mathbf{i})}_{h_\mfi,(D,\BQ_f)}(x)-\widetilde{f}^{(\mathbf{i})}_{h_\mfi,D}(x)\Big|\leq
\sup_{D\in\cQ_{\delta_n}}\Big|\widetilde{f}^{(\mathbf{i})}_{h_\mfi,(D,\BQ_f)}(x)-\varkappa_{h_{\mathbf{i}}}(D,\BQ_f,x)\Big|
\nonumber\\
\quad+\sup_{D\in\cQ_{\delta_n}}\Big|\widetilde{f}^{(\mathbf{i})}_{h_\mfi,D}(x)-\varkappa_{h_{\mathbf{i}}}(D,D,x)\Big|+
\sup_{D\in\cQ_{\delta_n}}\Big|\varkappa_{h_{\mathbf{i}}}(D,\BQ_f,x)-\varkappa_{h_{\mathbf{i}}}(D,D,x)\Big|
\\
\hskip-1.5cm\leq 2\chi_{\mathbf{i}}(f,x)+2\mathbf{C}L^2h_{\mathbf{i}}^{\beta}+\sup_{D\in\cQ_{\delta_n}}
\Big|\varkappa_{h_{\mathbf{i}}}(D,\BQ_f,x)-\varkappa_{h_{\mathbf{i}}}(D,D,x)\Big|
\\
\hskip-0.7cm\leq2\chi_{\mathbf{i}}(f,x)+2\mathbf{C}L^2h_{\mathbf{i}}^{\beta}+4CL^2\big(\cK,\beta,\sqrt{2}\big)h_{\mathbf{i}}^{\beta}
=2\chi_{\mathbf{i}}(f,x)+\mathbf{B}L^2h_{\mathbf{i}}^{\beta}.
\end{gather*}
To get the second inequality we have used Lemma \ref{lem:new1} while the last equality follows from the definition of $\mathbf{C}$.
Noting that the definition of $\widehat{Q}^{(\mathbf{i})}$ implies the inclusion $\bar{\cZ}^{(\mathbf{i})}\subseteq \{R^{(\mathbf{i})}_{n}(\BQ_f)\neq 0\}$
we obtain
\begin{eqnarray}
\label{eq20:proof}
\bar{\cZ}^{(\mathbf{i})}\subseteq \{R^{(\mathbf{i})}_{n}(\BQ_f)\neq 0\}&=&\Big\{\sup_{D\in\cQ_{\delta_n}}
\Big|\widetilde{f}^{(\mathbf{i})}_{h_\mfi,(D,\BQ_f)}(x)-\widetilde{f}^{(\mathbf{i})}_{h_\mfi,D}(x)\Big|>\mathbf{B}L^2h_{\mathbf{i}}^{\beta}\Big\}
\nonumber\\
&\subseteq&\big\{\chi_{\mathbf{i}}(f,x)\neq 0\big\}.
\end{eqnarray}
Denoting by $\omega_0=\ln(n)$, $\mathbf{n}_{\mathbf{0}}=\lfloor n/4\rfloor$ and

\vskip0.1cm

\centerline{$
\displaystyle{\mathbf{e}_{\mathbf{i}}=\big(L^{\frac{2}{\beta}}\omega_{\mathbf{i}}/\bn\big)^{-\frac{p\beta}{2\beta+1}}\sup_{\cF(\beta,L,\cQ_{\delta_n})}
\bE_f\big|\breve{f}^{(\mathbf{i})}(x)-f(x)\big|^p}.
$}

\vskip0.1cm

\noindent we deduce from    (\ref{eq19:proof}), (\ref{eq20:proof}) and the first assertion of Proposition \ref{prop4}

\vskip0.15cm

\centerline{$
\displaystyle{\mathbf{e}_{\mathbf{i}}\leq\mathbf{R}+\mathbf{P}\mathbf{e}_{\mathbf{i}-1}, \quad\forall \mathbf{i}=1,\ldots,\mathbf{i}^*,\;\forall n\geq 3}.
$}

\vskip0.15cm

\noindent
It yields together  with (\ref{eq18:proof}) since $\mathbf{P}<1$ for all $n\geq 3$

\vskip0.15cm

\centerline{$
\displaystyle{\mathbf{e}_{\mathbf{i}^*}\leq \mathbf{P}^{\mathbf{i}^*}\mathbf{e}_0+\mathbf{R}(1-\mathbf{P})^{-1}\leq \mathbf{e}_0 + 3^{p}\big\{(2\mathbf{B})^p+2^p\mathbf{E}+C_1 \big\}(1-\mathbf{P})^{-1}}.
$}

\vskip0.15cm

\noindent Since $\breve{f}^{(\mathbf{0})}(x)=\widehat{f}(x)$ we deduce from Theorem \ref{th2} that

\vskip0.1cm

\centerline{$
\displaystyle{\limsup_{n\to\infty}\mathbf{e}_0<\infty,}
$}

\vskip0.1cm

\noindent
that  completes the proof of the theorem.
\epr

\section{Proofs of Lemmas \ref{lem1}-\ref{lem:new2} and Proposition \ref{prop2}}

The proofs of Lemmas \ref{lem1}-\ref{lem:new2} are  based on the following result proved in the end of this section.

\begin{lemma}
\label{lem0}
For any $g\in\cG(\beta,L)$ and any $2\times 2$ matrix $\Psi=(\psi^T_1,\psi_2^T)$

\vskip0.1cm

\centerline{$
\displaystyle{\sup_{y\in\bR^2}\bigg|\int_{\bR^2}K(t)g\big(y+\Psi th\big)\rd t-g(y)\bigg|\rd \underline{u}\leq C(\cK,\beta,\psi^*)L^2h^{\beta},\quad\forall h>0,}
$}

\vskip0.1cm

\noindent
where $\psi^*=\|\psi_1\|\vee\|\psi_2\|$.

\end{lemma}

\noindent{\it Proof of Lemma \ref{lem1}.}
$\mathbf{1^0}.\;$  We obviously have
\begin{eqnarray*}
\cE_{h}(f,D,x)&:=&\bE_{f}\big[\widetilde{f}_{h,d}(x)\big]\bE_{f}\big[\widetilde{f}_{h,d_\perp}(x)\big]
\\
&=&\hskip-0.3cm\left[\int_{\bR^2}\cK_h\big(d^T(u-x)\big) f\big(u\big)\rd u\right]\left[\int_{\bR^2}\cK_h\big(d_\perp^T(u-x)\big) f\big(u\big)\rd u\right]
\\*[2mm]
&=:&
\cE^\prime_{h}(f,D,x)\cE^{\prime\prime}_{h}(f,D,x).
\end{eqnarray*}
Since $f\big(u\big)=\Bg_f\big(\BQ_f^Tu\big), u\in\bR^2$, denoting for brevity $\BQ_f=(\bq,\bq_\perp)$ and by $\Bg_i, i=1,2$, the marginals of
$\Bg_f$, we get
\begin{eqnarray*}
\label{eq2}
&&\cE^\prime_{h}(f,D,x)=\int_{\bR^2}\cK_h\big(u_1\big) \Bg_1\big(\bq^{T}x+\bq^{T}du_1+\bq^Td_\perp u_2\big)\times
\\
&&
\qquad\qquad\qquad\qquad\qquad\quad\;\; \Bg_2\big(\bq_\perp^{T}x+\bq_\perp^{T}du_1+\bq_\perp^Td_\perp u_2\big)\rd u_1\rd u_2
\nonumber\\
&&=\int_{\bR^2}\cK(s_1)  \Bg_1\big(\bq^{T}x+\bp_2s_1h+\bp_1 s_3\big) \Bg_2\big(q_\perp^{T}x-\bp_1s_1h+\bp_2 s_3\big)\rd s_1\rd s_3,
\\*[2mm]
&&\cE^{\prime\prime}_{h}(f,D,x)=
\int_{\bR^2}\cK_h\big(u_2\big) \Bg_1\big(\bq^{T}x+\bq^{T}du_1+\bq^Td_\perp u_2\big)
\times
\\
&&
\qquad\qquad\qquad\qquad\qquad\quad\;\;
\Bg_2\big(\bq_\perp^{T}x+\bq_\perp^{T}du_1+\bq_\perp^Td_\perp u_2\big)\rd u_1\rd u_2
\nonumber\\
&&=\int_{\bR^2}\cK(s_2)  \Bg_1\big(\bq^{T}x+\bp_2s_4+\bp_1 s_2h\big) \Bg_2\big(\bq_\perp^{T}x-\bp_1s_4+\bp_2 s_2h\big)\rd s_2\rd s_4.
\nonumber
\end{eqnarray*}
Thus we obtain that
\begin{eqnarray*}
&&\cE_{h}(f,D,x)
\\
&&=\int_{\bR^4}K(\overline{s})  \Bg_f\big(\BQ_f^{T}x+\bp_2\overline{s}h+\bp_1\Gamma\underline{s}\big)
\Bg_f\big(\BQ_f^{T}x+\bp_1\Gamma\Omega \overline{s}h+\bp_2\Omega\underline{s}\big) \rd s.
\end{eqnarray*}
If $D=\BQ_f$ that implies $\bp_1=0$ and $\bp_2=1$ we get
\begin{eqnarray*}
\cE_{h}(f,D,x)
&=&\int_{\bR^4}K(\overline{s})  \Bg_f\big(\BQ_f^{T}x+\overline{s}h\big)
\Bg_f\big(\BQ_f^{T}x+\Omega\underline{s}\big) \rd s
\\
&=&
\int_{\bR^2}K(\overline{s})  \Bg_f\big(\BQ_f^{T}x+\overline{s}h\big)\rd\overline{s},
\end{eqnarray*}
since $\Bg_f$ is a probability density. The assertion of the lemma in this case follows from Lemma \ref{lem0}. If $D\neq\BQ_f$ ($\bp_1\neq 0$),
making the change of variables
$
\BQ_f^{T}x+\bp_1\Gamma\underline{s}=\bp_1\Gamma\underline{t}
$
and noting that $\Gamma^{-1}=\Gamma$ we come to
\vskip-0.6cm
\begin{gather*}
\label{eq400000}
\hskip-10.2cm \cE_{h}(f,D,x)
\\
=\int_{\bR^4}K(\overline{t}) \Bg_f\big(\bp_2\overline{t}h+\bp_1\Gamma\underline{t}\big) \Bg_f\big([I-\bp_2\bp_1^{-1}\Omega\Gamma]\BQ_f^{T}x+\bp_1\Gamma\Omega \overline{t}h+\bp_2\Omega\underline{t}\big)\rd t.
\end{gather*}
Noting that
$
I-\bp_2\bp_1^{-1}\Omega\Gamma=
\bp_1^{-1}D^T\BQ_f\Gamma\Omega,
$
we get

\vskip0.1cm

\centerline{$
[I-\bp_2\bp_1^{-1}\Omega\Gamma]\BQ_f^{T}x
=\bp_1^{-1}
\left(
\begin{array}{cccc}
d_\perp^Tx
\\*[3mm]
- d^Tx
\end{array}\right).
$}

\vskip0.1cm

\noindent Thus we have
\vskip-0.7cm
\begin{eqnarray*}
&&\Bg_f\big([I-\bp_2\bp_1^{-1}\Omega\Gamma]\BQ_f^{T}x+\bp_1\Gamma\Omega \overline{t}h+\bp_2\Omega\underline{t}\big)
\\
&&=\Bg_1\big(\bp_1^{-1}d_\perp^Tx+\bp_1
t_2h+\bp_2t_4\big)\Bg_2\big(-\bp_1^{-1}d^Tx-\bp_1t_1h+\bp_2t_3\big)
\end{eqnarray*}
\vskip-0.1cm
\noindent and, since $\Bg_2$ is symmetric
\vskip-0.6cm
\begin{eqnarray*}
&&\Bg_f\big([I-\bp_2\bp_1^{-1}\Omega\Gamma]\BQ_f^{T}x+\bp_1\Gamma\Omega \overline{t}h+\bp_2\Omega\underline{t}\big)
\\
&&=\Bg_1\big(\bp_1^{-1}d_\perp^Tx+\bp_1
t_2h+\bp_2t_4\big)\Bg_2\big(\bp_1^{-1}d^Tx+\bp_1t_1h-\bp_2t_3\big).
\end{eqnarray*}
\vskip-0.1cm
\noindent Noting that
$
\left(
\begin{array}{cccc}
d_\perp\;
d
\end{array}\right)^T=D\Omega,
$
 we obtain finally
 \vskip-0.4cm
\begin{eqnarray}
\label{eq40000000}
&&\cE_{h}(f,D,x)
\nonumber\\
&&=\int_{\bR^4}K(\overline{t}) \Bg_f\big(\bp_2\overline{t}h+\bp_1\Gamma\underline{t}\big) \Bg_f\big(\bp_1^{-1}D\Omega x+\bp_1\Omega \overline{t}h+\bp_2\Omega\Gamma\underline{t}\big)\rd t.
\end{eqnarray}
Consider now two cases.

$\mathbf{2^0a}.\;$ If $|\bp_2|\geq |\bp_1|$  using  $\Omega D\Omega=D$, $\Omega^2=I, \Gamma^2=I$ and   making the change of variables $\overline{t}=\overline{v}$,

\vskip0.15cm

\centerline{$
\bp_1^{-1}D\Omega x+\bp_1\Omega \overline{t}h+\bp_2\Omega\Gamma\underline{t}=\underline{v}\;\Rightarrow\;\underline{t}=\bp_2^{-1}\Gamma\Omega\underline{v}-\bp_2^{-1}\big[\bp_1^{-1}\Gamma D x+\bp_1\Gamma \overline{v}h\big]
$}

\vskip0.15cm

\noindent we obtain (remind that $\Bg_f$ is a symmetric function)
\vskip-0.65cm
\begin{eqnarray*}
\label{eq40}
&&\cE_{h}(f,D,x)
\\
&&=\bp_2^{-2}\int_{\bR^4}\Bg_f\big(\underline{v}\big)K(\overline{v}) \Bg_f\big(\bp_2^{-1} D x-\big[\bp_2-\bp^2_1\bp_2^{-1}\big]\overline{v}h-\bp_1\bp_2^{-1}\Omega\underline{v}\big) \rd v.
\end{eqnarray*}
\vskip-0.15cm
\noindent Hence, taking into account that $\Bg_f$ is a probability density we deduce from Lemma \ref{lem0} that
\vskip-0.55cm
\begin{gather}
\label{eq400}
\qquad\Big|\cE_{h}(f,D,x)
-\bp_2^{-2}\int_{\bR^2}\Bg_f\big(\underline{v}\big) \Bg_f\big(\bp_2^{-1} D x-\bp_1\bp_2^{-1}\Omega\underline{v}\big) \rd \underline{v}\Big|
\\
=\bp_2^{-2}\bigg|\int_{\bR^2}\Bg_f\big(\underline{v}\big)\bigg[\int_{\bR^2}K(\overline{v}) \Bg_f\big(\bp_2^{-1} D x-\big[\bp_2-\bp^2_1\bp_2^{-1}\big]\overline{v}h-\bp_1\bp_2^{-1}\Omega\underline{v}\big) \rd \overline{v}
\nonumber\\
\qquad\qquad\qquad\quad
-\int_{\bR^2} \Bg_f\big(\bp_2^{-1} D x-\bp_1\bp_2^{-1}\Omega\underline{v}\big)\bigg] \rd \underline{v}\bigg|
\nonumber\\
\leq 2 \int_{\bR^2}\Bg_f(\underline{v})\sup_{y\in\bR^2}\bigg|\int_{\bR^2}K(\overline{v})\Bg_f\big(y-\big[\bp_2-\bp^2_1\bp_2^{-1}\big] \overline{v}h\big)\rd \overline{v}-\Bg_f(y)\bigg|\rd \underline{v}
\nonumber\\*[2mm]
\hskip-7.8cm\leq 2C(\cK,\beta,1)L^2h^{\beta}.
\nonumber
\end{gather}
\vskip-0.1cm
\noindent Here we have also used that $\bp_1^2+\bp_2^2=1$ and therefore $(\bp_1\vee \bp_2)^2\geq 1/2$.

\vskip0.1cm

$\mathbf{2^0b}$ If $|\bp_2|< |\bp_1|$ making the change of variables $\overline{t}=\overline{v}$ and $\bp_2\overline{t}h+\bp_1\Gamma\underline{t}=\underline{v}$, we obtain
\begin{eqnarray*}
&&\cE_{h}(f,D,x)
\\
&&=\bp_1^{-2}\int_{\bR^4}\Bg_f\big(\underline{v}\big)K(\overline{v}) \Bg_f\big(\bp_1^{-1}D\Omega x+[\bp_1-\bp_2^2\bp_1^{-1}]\Omega \overline{v}h+\bp_2\bp_1^{-1}\Omega\underline{v}\big) \rd v.
\end{eqnarray*}
We deduce from Lemma \ref{lem0} similarly to (\ref{eq400})
\begin{eqnarray}
\label{eq401}
&&\Big|\cE_{h}(f,D,x)
-\bp_1^{-2}\int_{\bR^2}\Bg_f\big(\underline{v}\big) \Bg_f\big(\bp_1^{-1}D\Omega x+\bp_2\bp_1^{-1}\Omega\underline{v}\big) \rd v\Big|
\nonumber\\
&&\leq 2C(\cK,\beta,1)L^2h^{\beta}.
\end{eqnarray}
It is worth noting that (\ref{eq400}) and (\ref{eq401}) can be written in a unified way
\begin{eqnarray*}
\big|\cE_{h}(f,D,x)-\cE_{0}(f,D,x)\big|
\leq 2C(\cK,\beta,1)L^2h^{\beta}.
\end{eqnarray*}
Thus, remarking that
$
\tau_f\big(D,\BQ_f\big)=\cE_{0}(f,D,x)
$
we come to the assertion of the lemma.
\epr

\noindent {\it Proof of Lemma \ref{lem:new1}.} Since by definition  $\widetilde{f}_{h,(\BQ_f,\BQ_f)}(x)=\widetilde{f}_{h,\BQ_f}(x)$
it suffices to prove the lemma for any $D\neq \BQ_f$.
We obviously have
\begin{eqnarray*}
&&E_{h}(f,D,x):=\bE_{f}\big[\widetilde{f}_{h,(D,\BQ_f)}(x)\big]
\\
&=&\int_{\bR^4}K_{h}(z)\Bg_f\big(\BQ_f^T\overline{y}+\bp_1\BQ_f^T\Gamma\Omega \overline{z}+\bp_2\BQ_f^T\Omega\underline{z} \big)
\Bg_f\big(\bp_2\BQ_f^T \overline{z}+\bp_1\BQ_f^T\Gamma\underline{z} \big)\rd z.
\end{eqnarray*}
Noting that

\vskip0.05cm

\centerline{$
\BQ_f^T\Omega \BQ_f^T=\Omega,\quad \BQ_f^T\Gamma \BQ_f^T=\Gamma
$}

\vskip0.15cm

\noindent and putting  $\underline{z}=\BQ_f^T\underline{u}, \overline{z}=\overline{u}h$ we get
\vskip-0.7cm
\begin{eqnarray*}
&&E_{h}(f,D,x)
\\
&&=\int _{\bR^4}K\big(\overline{u}\big)\Bg_f\big(\BQ_f^T\overline{y}+\bp_1\BQ_f^T\Gamma\Omega \overline{u}h+\bp_2\Omega\underline{u} \big)
\Bg_f\big(\bp_2\BQ_f^T \overline{u}h+\bp_1\Gamma\underline{u} \big)\rd u.
\end{eqnarray*}
\vskip-0.2cm
Consider now two cases.

$\mathbf{1^0a}.\;$ If $|\bp_2|\geq |\bp_1|$  using   $\Omega^2=I$, $\Gamma\Omega \BQ_f^T\Gamma\Omega=-\BQ_f^T$, $\bp_1^2+\bp_2^2=1$ and   making the change of variables $\overline{u}=\overline{v}$,
\vskip-0.5cm
$$
\BQ_f^T\overline{y}+\bp_1\BQ_f^T\Gamma\Omega \overline{u}h+\bp_2\Omega\underline{u}=\underline{v}\;\Rightarrow\;\underline{u}=\bp_2^{-1}\Omega\underline{v}-\bp_2^{-1}\Omega \BQ_f^T y-\bp_1\bp_2^{-1}\Omega \BQ_f^T\Gamma\Omega  \overline{v}h
$$
\vskip-0.2cm
\noindent we obtain
\vskip-0.7cm
\begin{gather*}
\hskip-10.2cm E_{h}(f,D,x)
\\
=\bp_2^{-2}\int _{\bR^4}\Bg_f\big(\underline{v} \big)K\big(\overline{v}\big)
\Bg_f\big(-\bp_1\bp_2^{-1}\BQ_f^T\overline{y}\Gamma\Omega + \bp_1\bp_2^{-1}\Gamma\Omega\underline{v}+\bp_2^{-1}\BQ_f^T\overline{v}h \big)\rd v.
\end{gather*}
\vskip-0.2cm
\noindent Applying Lemma \ref{lem0} we obtain similarly to (\ref{eq400})
\vskip-0.7cm
\begin{eqnarray}
\label{eq5}
&&\bigg| E_{h}(f,D,x)-\bp_2^{-2}\int _{\bR^2}\Bg_f\big(\underline{v} \big)
\Bg_f\big(-\bp_1\bp_2^{-1}\BQ_f^T\overline{y}\Gamma\Omega + \bp_1\bp_2^{-1}\Gamma\Omega\underline{v}\big)\rd \underline{v}\bigg|
\nonumber\\
&&\leq 2C\big(\cK,\beta,\sqrt{2}\big)L^2h^{\beta}.
\end{eqnarray}
\quad$\mathbf{1^0b}.\;$ If $|\bp_1|> |\bp_2|$  using   $\Gamma^2=I$, $\BQ_f^T\Gamma\Omega=-\Omega\Gamma \BQ_f^T$, $\bp_1^2+\bp_2^2=1$ and   making the change of variables $\overline{u}=\overline{v}$,

\vskip0.15cm

\centerline{$
\bp_2\BQ_f^T \overline{u}h+\bp_1\Gamma\underline{u} =\underline{v}\;\;\Rightarrow\;\;\underline{u}=\bp_1^{-1}\Gamma\underline{v}-\bp_2\bp_1^{-1}\Gamma \BQ_f^T \overline{v}h
$}

\vskip0.15cm

\noindent we obtain
\vskip-0.7cm
\begin{eqnarray*}
&&E_{h}(f,D,x)
\\
&&=\bp_1^{-2}\int _{\bR^4}\Bg_f\big(\underline{v} \big)K\big(\overline{v}\big)
\Bg_f\big(\BQ_f^T\overline{y} + \bp_2\bp_1^{-1}\Omega\Gamma\underline{v}+\bp_1^{-1}\BQ_f^T\Gamma\Omega\overline{v}h \big)\rd v.
\end{eqnarray*}
\vskip-0.3cm
\noindent The application of Lemma \ref{lem0} yields
\vskip-0.7cm
\begin{eqnarray}
\label{eq6}
&&\Big| E_{h}(f,D,x)-\bp_1^{-2}\int _{\bR^2}\Bg_f\big(\underline{v} \big)
\Bg_f\big(\BQ_f^T\overline{y} + \bp_2\bp_1^{-1}\Omega\Gamma\underline{v}\big)\rd \underline{v}\Big|
\nonumber\\
&&
\leq 2C\big(\cK,\beta,1\big)L^2h^{\beta}\leq 2C\big(\cK,\beta,\sqrt{2}\big)L^2h^{\beta} .
\end{eqnarray}
\vskip-0.2cm
\noindent Note that (\ref{eq5}) and (\ref{eq6}) can be written as

\vskip0.15cm

\centerline{$
\big|E_{h}(f,D,x)-E_{0}(f,D,x)\big|
\leq 2C\big(\cK,\beta,\sqrt{2}\big)L^2h^{\beta},
$}

\vskip0.15cm

\noindent Note also that $\BQ_f^T\overline{y}=\bp^{-1}_1D\Omega x$ and
\vskip-0.4cm
\begin{eqnarray*}
E_{0}(f,D,x)&=&\int _{\bR^2}\Bg_f\big(\bp^{-1}_1D\Omega x+\bp_2\Omega\underline{u} \big)
\Bg_f\big(\bp_1\Gamma\underline{u} \big)\rd \underline{u}
\\
&=&\int_{\bR^2} \Bg_f\big(\bp_1^{-1}D\Omega x+\bp_2\Omega\Gamma u\big)\Bg_f\big(\bp_1\Gamma u\big)\rd u=\tau_f(D,\BQ_f).
\end{eqnarray*}
\vskip-0.2cm
\noindent To get the penultimate equality we used the change of variables  $u_1=\underline{v}_1$, $u_2=-\underline{v}_2$ and the symmetry of $\Bg_{f,2}$ which implies $g\big(p_1 v \big)=g\big(p_1\Gamma v \big)$. Hence,

\vskip0.1cm

\centerline{$\Big|\bE_{f}\big[\widetilde{f}_{h,(D,\BQ_f)}(x)\big]-\tau_f(D,\BQ_f)\Big|
\leq 2C\big(\cK,\beta,\sqrt{2}\big)L^2h^{\beta}
$}

\vskip0.1cm

\noindent that implies together with Lemma \ref{lem1} the assertion of the lemma.
\epr

\noindent {\it Proof of Lemma \ref{lem:new2}.}
 As it was  mentioned in Remark \ref{rem1}  $p_1(D,Q)=-p_1(Q,D)$ and  $p_2(D,Q)=p_2(Q,D)$. Moreover  $DQ=QD$ for any $D,Q\in\mQ$.
Hence
\begin{eqnarray*}
\widetilde{f}_{h,(Q,D)}(x)&=&\frac{1}{n(n-1)}\sum_{k,l=1,k\neq l}^nK_{h}\big(-p_1\Omega\Gamma X_k+p_2X_l-\Gamma\Omega QD\Omega x\big)
\\
&=&\frac{1}{n(n-1)}\sum_{k,l=1,k\neq l}^nK_{h}\big(p_1\Omega\Gamma (-X_k)+p_2X_l-\Gamma\Omega QD\Omega x\big).
\end{eqnarray*}
Remind that the density of $X_k$ is $\Bg_f\big(\BQ_f^Tv\big)$ and therefore the law of $-X_k$ coincides with whose of $X_k$ because $\Bg_f$ is symmetric.
Finally since $X_k$ and $X_l$ are independent for all $k\neq l$ for any $D,Q\in\mQ$
we conclude that

\vskip0.1cm

\centerline{$
K_{h}\big(-p_1\Omega\Gamma X_k+p_2X_l-\Gamma\Omega QD\Omega x\big)\stackrel{\text{law}}{=}K_{h}\big(p_1\Omega\Gamma X_k+p_2X_l-\Gamma\Omega QD\Omega x\big).
$}

\vskip0.15cm

\noindent It implies in particular the assertion of the lemma.
\epr

\noindent{\it Proof of Proposition \ref{prop2}.}  Denoting $(\xi_{1,i},\xi_{2,i})^T=\BQ^T_f(X_i-x), i=1,\ldots,n,$ we remark that

\vskip0.2cm

\centerline{$
\widetilde{f}^{(\cI)}_{\mh,\BQ_f}(x)=\Big[|\cI|^{-1}\sum_{k\in\cI}\cK_{\blh}\big(\xi_{1,k}\big)\Big]
\Big[|\cI|^{-1}\sum_{k\in\cI}\cK_{\blh}\big(\xi_{2,k}\big)\Big]=:\Upsilon_1(\blh)\Upsilon_2(\blh).
$}

\vskip0.2cm

\noindent Note that $\xi_{1,i},\xi_{2,i}, i=1,\ldots,n$ are independent with the densities given by
$\Bg_{1}\big(\bullet+\bq^Tx\big)$ and $\Bg_{2}\big(\bullet+\bq_\perp^Tx\big)$ respectively. We obviously have
\begin{eqnarray*}
\widetilde{f}^{(\cI)}_{\blh,\BQ_f}(x)-f(x)&=&\big[\Upsilon_1(\blh)-\bE_{\Bg_f}\big\{\Upsilon_1(\blh)\big\}\big]
\big[\Upsilon_2(\mh)-\bE_{\Bg_f}\big\{\Upsilon_2(\blh)\big\}\big]
\nonumber\\
&&\quad+\bE_{\Bg_f}\big\{\Upsilon_1(\blh)\big\}\big[\Upsilon_2(\blh)-\bE_{\Bg_f}\big\{\Upsilon_2(\blh)\big\}\big]
\\
&&\quad+\bE_{\Bg_f}\big\{\Upsilon_2(\blh)\big\}\big[\Upsilon_1(\blh)-\bE_{\Bg_f}\big\{\Upsilon_1(\blh)\big\}\big]
\nonumber
\\
&&\quad+\bE_{\Bg_f}\big\{\Upsilon_1(\blh)\big\}\bE_{\Bg_f}\big\{\Upsilon_2(\blh)\big\}-f(x).
\nonumber
\end{eqnarray*}
Here $\bE_{\Bg_f}$ is  the expectation w.r.t  the law of $\xi_1,\ldots,\xi_n$. In view of Lemma \ref{lem1}
\vskip-0.5cm
\begin{gather}
\Big|\bE_{\Bg_f}\big\{\Upsilon_1(\blh)\big\}\bE_{\Bg_f}\big\{\Upsilon_1(\blh)\big\}-f(x)\Big|=\Big|\bE_{\Bg_f}\big\{\cK_{\blh}\big(\xi_{1,1}\big)\big\}
\bE_{\Bg_f}\big\{\cK_{\blh}\big(\xi_{2,1}\big)\big\}-f(x)\Big|
\nonumber\\
\label{eq11:proof}
\quad\leq 2C(\cK,\beta,1)L^2 \blh^\beta=2C(\cK,\beta,1)L^{2}\big(\mu/|\cI|\big)^{\frac{\beta}{2\beta+1}}.
\end{gather}
Since $\Bg_f\in\cG(\beta,L)$ it implies  $\Bg_{1}$, $\Bg_{2}$ are uniformly bounded by $L$. Hence
\begin{eqnarray*}
&&\big|\bE_{\Bg_f}\big\{\Upsilon_j(\blh)\big\}\big|\leq L\|\cK\|_1,\qquad \bV_{\Bg_f}\big\{|\cI|\Upsilon_j(\blh)\big\}\leq L\|\cK\|^2_2|\cI|\blh^{-1},\; j=1,2;
\\*[1mm]
&&\bE_{\Bg_f}\big\{\big|\cK_{\blh}\big(\xi_{1,k}\big|^p\big\}\leq L\|\cK\|_p^p\blh^{1-p},\quad \bE_{\Bg_f}\big\{\big|\cK_{\blh}\big(\xi_{2,k}\big|^p\big\}\leq L\|\cK\|_p^p\blh^{1-p}.
\end{eqnarray*}
Applying the Rosenthal inequality (if $p>2)$  to  $|\cI|[\Upsilon_j(\mh)-\bE_{\Bg_f}\big\{\Upsilon_j(\mh)\big\}], j=1,2,$ which is a sum of i.i.d bounded and centered random variables or computing its variance (if $1\leq p\leq 2$) we assert that there exists $C>0$ completely determined by $p$ and $\cK$ such that
for any $n\geq 1$, $\cI\in\mI_n$ and $\mu\geq 1$
\begin{eqnarray}
\label{eq12:proof}
&&\bE_{\Bg_f}\big\{\big|\Upsilon_j(\blh)\big|^p\big\}\leq C(L^{p/2}+L)\big(\mu/|\cI|\big)^{\frac{\beta}{2\beta+1}},\;j=1,2.
\end{eqnarray}
The assertion of the proposition follows now from (\ref{eq11:proof}) and (\ref{eq12:proof}).
\epr

\paragraph{Proof of Lemma \ref{lem0}} Remind that for any function $w\in\bH(\beta,L)$
\vskip-0.6cm
\begin{eqnarray}
\label{eq42}
\sup_{z,\mz\in\bR}|z-\mz|^{-\beta}\bigg|\sum_{j=0}^m\frac{w^{(j)}(\mz)(z-\mz)^{j}}{j!}-w(z)\bigg|\leq L.
\end{eqnarray}
\vskip-0.2cm
\noindent We deduce from (\ref{eq42}) for any $t\in\bR^2$
\begin{eqnarray*}
\label{eq3-new}
\bigg|g_1(y_1+h\psi^T_1t)-\sum_{j=0}^m\frac{g_1^{(j)}(y_1)h^j(\psi^T_1t)^{j}}{j!}\bigg|&\leq& Lh^\beta|\psi^T_1t|^{\beta}\leq \psi^*Lh^\beta\|t\|^\beta ;
\nonumber\\
\bigg|g_2(y_2+h\psi^T_2t)-\sum_{j=0}^m\frac{g_2^{(j)}(y_2)h^j(\psi^T_2t)^{j}}{j!}\bigg|&\leq& Lh^\beta|\psi^T_2t|^{\beta}\leq \psi^*Lh^\beta\|t\|^\beta,
\end{eqnarray*}
where $\|\cdot\|$ is used for the euclidian norm. Setting

\vskip0.0cm

\centerline{$
\displaystyle{\mathbf{P}_{g,\Psi,y}(t)=\sum_{j,s=0}^m\frac{g_1^{(j)}(y_1)g_2^{(s)}(y_2)h^{j+s}(\psi^T_1t)^{j}(\psi^T_2t)^{s}}{j!s!}}
$}

\vskip0.1cm

\noindent and recalling that $\|g_i\|_\infty\leq L, i=1,2,$ we obviously have
\begin{eqnarray}
\label{eq44}
\big|g\big(y+\Psi th\big)-\mathbf{P}_{g,\Psi,y}(t)\big|\leq 2\psi^*L^2h^\beta\|t\|^\beta+(\psi^*L)^2h^{2\beta}\|t\|^{2\beta}.
\end{eqnarray}
It remains to note that $\mathbf{P}_{g,\Psi,y}(t)$ can be rewritten as

\vskip0.1cm

\centerline{$
\displaystyle{
\mathbf{P}_{g,\Psi,y}(t)=\sum_{i,l=0}^{2m}a_{i,l}t_1^it_2^l,\quad a_{0,0}=g_1(y_1)g_2(y_2)=g(y)},
$}

\vskip0.1cm

\noindent
and, therefore, in view of Assumption \ref{ass:kernel}

\vskip0.1cm

\centerline{$
\displaystyle{
\int_{\bR^2}K(t)\mathbf{P}_{g,\Psi,y}(t)\rd t=g(y)}.
$}

\vskip0.1cm

\noindent
This together with (\ref{eq44}) allows as to assert that

\vskip0.1cm

\centerline{$
\displaystyle{
\bigg|\int_{\bR^2}K(t)g\big(y+\Psi th\big)\rd t-g(y)\bigg|\leq C(\cK,\beta,\psi^*)L^2h^{\beta},\quad\forall h>0.}
$}

\vskip0.1cm

\noindent
Lemma is proved.
\epr

\section{Proofs of Propositions \ref{prop1}-\ref{prop4}} Set $\beta_1,\beta_2>0$, $L\geq 1$ and let $f\in\cF(\beta,L,\cQ_{\delta_n})$, $\beta\in\{\beta_1,\beta_2\}$, be fixed. We divide these proofs into three steps.

\smallskip

\noindent\textit{First step: upper bounds for sums of independent variables.}

\vskip0.1cm

\noindent For any $(h,D)\in\bR_+^*\times\cQ_{\delta_n}$ and any $b\in\{d,d_{\perp}\}$ set
\begin{eqnarray*}
&& \xi_{h,D}(x):=\widetilde{f}_{h,D}(x)-\bE_f\{\widetilde{f}_{(h,d)}(x)\}\bE_f\{\widetilde{f}_{(h,d_\perp)}(x)\},
\\
&& \xi_{(h,b)}(x):=\widetilde{f}_{(h,b)}(x)-\bE_f\{\widetilde{f}_{(h,b)}(x)\},
\\
&& G_{h,b}(x):=1\vee\bE_f\left\{\left|\cK_h\left(b^TX_1-b^Tx\right)\right|\right\},
\\
&&\widetilde{G}_{h,b}(x):=1\vee\Big[\frac{1}{n}\sum_{k=1}^n\left|\cK_{h}\big(b^TX_k-b^Tx\big)\right|\Big].
\end{eqnarray*}

\vskip0.1cm

\noindent Note first that, since $|\bq^Tb|^2+|\bq^Tb_\perp|^2=1$ and $L\wedge\|\cK\|_1\geq 1$,
\begin{eqnarray*}
&&G_{h,b}(x)= 1\vee\int_{\bR^2}\left|\cK(v_1)\right|\Bg_1\left(\bq^Tx+h\bq^Tbv_1+\bq^Tb_\perp v_2\right)\times
\\
&&\qquad\qquad\qquad\qquad\qquad\quad\Bg_2\left(\bq_\perp^Tx+h\bq_\perp^Tbv_1+\bq_\perp^Tb_\perp v_2\right)\rd v\leq\sqrt{2}\|\cK\|_1L.
\end{eqnarray*}

\vskip0.1cm

\noindent  For any $q\geq 1$ and any $\epsilon >0$ put $\displaystyle{\lambda_q^{(1)}(\epsilon)=\left[\sqrt{2}+\sqrt{5q\epsilon^{-1}}\right](\|\cK\|_\infty\vee 1)}$. Consider finally a real number $\displaystyle{\alpha_n\geq 1\vee\mn(\cQ_{\delta_n})}$. In the sequel $\alpha_n$ and $\epsilon$ will be fixed and properly chosen.

\vskip0.1cm

\noindent Applying Bernstein inequality we obtain for any $q\geq 2$, any integer $n\geq3$, any  $z\in [0, 2q\alpha_n]$ and all real numbers $h$ satisfying $nh\geq\epsilon\alpha_n$
\begin{eqnarray}
\label{uperbound1}
\sup_{b\in\{d,d_\perp\}}\bP_f\bigg\{\sup_{D\in\cQ_{\delta_n}}\bigg[\left|\xi_{(h,b)}(x)\right|-\lambda_q^{(1)}(\epsilon)G_{h,b}(x)
\sqrt{\frac{0.5q\alpha_n+z}{nh}}\;\bigg]>0\bigg\}
\\
\quad\leq 2 e^{-z}.
\nonumber
\end{eqnarray}
By integration of the Bernstein inequality we get for any $q\geq 1$,  $n\geq3$, any  $t\in [0, 1.5q\alpha_n]$ and any real  $h$ satisfying $nh\geq\epsilon\alpha_n$
\begin{eqnarray}
\label{uperbound2}
\sup_{b\in\{d,d_\perp\}}\bE_f\bigg\{\sup_{D\in\cQ_{\delta_n}}\bigg[\left|\xi_{(h,b)}(x)\right|-
\lambda_q^{(1)}(\epsilon)G_{h,b}(x)\sqrt{\frac{q\alpha_n+t}{nh}}\;\bigg]_+\bigg\}^q
\\
\leq C_1^{(q)}(\cK)L^q\left[nh\right]^{-\frac{q}{2}} e^{-t},\nonumber
\end{eqnarray}
where $C_1^{(q)}(\cK,\epsilon)=2^{q+1}\big(1+\sqrt{\epsilon^{-1}}\big)^q\Gamma(q+1)\|\cK\|_1^q\left(\|\cK\|_\infty\vee 1\right)^q$ and $\Gamma$ is the Gamma function.
 Choose now $\alpha_n=\alpha\ln(n)$, $\epsilon=1$ and $\;t=1,5q\alpha_n$ and  introduce  $\gamma_q=\lambda_q^{(1)}(1)\sqrt{2.5q\alpha}$,

\vskip0.1cm

\centerline{$\mH_q:=\Big\{h\in\cH :\; nh\geq 10q\big[\lambda^{(1)}_q(1)\big]^2\alpha_n\Big\}$.}

\vskip0.1cm

\noindent Since $\text{card}(\mH_q)\leq\ln(n)$, we deduce from (\ref{uperbound2}) that for any $q\geq 1$ and $n\geq 3$
\begin{eqnarray*}
&& \sup_{b\in\{d,d_\perp\}}\bE_f\bigg\{\sup_{h\in\mH_q}\sup_{D\in\cQ_{\delta_n}}\left[\left|\xi_{(h,b)}(x)\right|-\gamma_qG_{h,b}(x)\sqrt{\ln(n)/nh}\;\right]_+
\bigg\}^q
\\
&&\leq C_1^{(q)}(\cK,1)L^q\left[\ln(n)/n\right]^{\frac{q}{2}}n^{-q}.
\end{eqnarray*}
 Additionally, using kernel $|\cK|$ instead of $\cK$ in the last inequality we get
\begin{eqnarray*}
&&\sup_{b\in\{d,d_\perp\}}\bE_f\bigg\{\sup_{h\in\mH_q}\sup_{D\in\cQ_{\delta_n}}\left[G_{h,b}(x)-2\widetilde{G}_{h,b}(x)\;\right]_+\bigg\}^q
\\
&&\leq 2^qC_1^{(q)}(\cK,1)L^q\left[\frac{\ln(n)}{n}\right]^{\frac{q}{2}}n^{-q},
\\
&&\bE_f\bigg\{\sup_{h\in\mH_q}\sup_{D\in\cQ_{\delta_n}}\sup_{b\in\{d,d_\perp\}}\left[\widetilde{G}_{h,b}(x)\right]\bigg\}^q
 \\
&& \leq 2^{q-1}\left[2C_1^{(q)}(\cK,1)+\left(\frac{3}{\sqrt{2}}\right)^q\|\cK\|_1^q\right]L^q.
\end{eqnarray*}

\noindent Noting that $\displaystyle{\widehat{\cU}_n=\sup_{h\in\mH}\sup_{D\in\cQ_{\delta_n}}\sup_{b\in\{d,d_\perp\}}\left[\widetilde{G}_{h,b}(x)\right]^2}$ and that $\mH\subset\mH_{2p}$ for $n$ large enough, we obtain the third assertion of Proposition \ref{prop4} with

\vskip0.1cm

\centerline{$\mathbf{C}_p(\cK)=2^{2p-1}\left[2C_1^{(2p)}(\cK,1)+(3/\sqrt{2})^{2p}\|\cK\|_1^{2p}\right]$.}

\vskip0.1cm

\noindent Since $\displaystyle{|\xi_{h,D}(x)|\leq|\xi_{(h,d)}(x)|\times\widetilde{G}_{h,d_\perp}(x)+|\xi_{(h,d_\perp)}(x)|\times G_{h,d}(x)}$, using Cauchy-Schwartz inequality we get for all $p\geq 1$ and all $n$ large enough
\begin{eqnarray}
\label{uperbound3}
&&\bE_f\bigg\{\sup_{h\in\mH}\sup_{D\in\cQ_{\delta_n}}\left[\left|\xi_{h,D}(x)\right|-\ma\widehat{\cU}_n\sqrt{\ln(n)/nh}\;\right]_+\bigg\}^p
\\
&&\leq CL^{2p}\left[\ln(n)/n\right]^{\frac{p}{2}}n^{-p},\quad C:=C(p,\cK)>0.\nonumber
\end{eqnarray}

\noindent Similarly, in view of (\ref{uperbound1}) with $q=2p$ and $z=4p\alpha_n$, one has for all $p\geq 1$ and all $n$ large enough
\begin{eqnarray}
\label{uperbound4}
 \sup_{b\in\{d,d_\perp\}}\bP_f\bigg\{\sup_{h\in\mH}\sup_{D\in\cQ_{\delta_n}}\left[\left|\xi_{(h,b)}(x)\right|-
\gamma_{2p}G_{h,b}(x)\sqrt{\ln(n)/nh}\;\right]>0\bigg\}\nonumber
\\
 \leq 2\ln(n)n^{-4p},\nonumber
\\*[1mm]
\bP_f\bigg\{\sup_{h\in\mH}\sup_{D\in\cQ_{\delta_n}}\left[\left|\xi_{h,D}(x)\right|-\ma\widehat{\cU}_n\sqrt{\ln(n)/nh}\;\right]>0\bigg\}\leq 10\ln(n)n^{-4p}.
\end{eqnarray}

\noindent\textit{Second step: upper bounds for $U$-Statistics of Order Two.}

\noindent For any $(D,Q)\in\cQ_{\delta_n}^2$, $Q\neq D$, and any $h>0$ put
\begin{eqnarray*}
\xi_{h,(D,Q)}(x) & := & \overline{f}_{h,(D,Q)}(x)-\bE_f\left[\overline{f}_{h,(D,Q)}(x)\right]\quad\text{and}
\\
\varphi(X_k,X_l) & := & \frac{1}{n(n-1)}K_{h}\big(p_1\Omega\Gamma X_k+p_2X_l-\Omega\Gamma QD\Omega x\big)
\end{eqnarray*}

\noindent Let's write $\displaystyle{\xi_{h,(D,Q)}(x) = \xi_{h,(D,Q)}^{(1)}(x)+\xi_{h,(D,Q)}^{(2)}(x)+\xi_{h,(D,Q)}^{(3)}(x)}$, where
\begin{eqnarray*}
 \xi_{h,(D,Q)}^{(1)}(x) & := & \sum_{k,l=1,\; k\neq l}^n\big(\varphi(X_k,X_l)-\bE_f\left[\varphi(X_k,X_l)|X_l\right]
\\
&& \qquad\quad-\bE_f\left[\varphi(X_k,X_l)|X_k\right]+\bE_f\left[\varphi(X_k,X_l)\right]\big),
\\
\xi_{h,(D,Q)}^{(2)}(x) &:= & \sum_{k,l=1,\; k\neq l}^n\left(\bE_f\left[\varphi(X_k,X_l)|X_l\right]-\bE_f\left[\varphi(X_k,X_l)\right]\right),
\\*[2mm]
\xi_{h,(D,Q)}^{(3)}(x) &:= &  \sum_{k,l=1,\; k\neq l}^n\left(\bE_f\left[\varphi(X_k,X_l)|X_k\right]-\bE_f\left[\varphi(X_k,X_l)\right]\right).
\end{eqnarray*}
Note that $\;\displaystyle{\xi_{h,(D,Q)}^{(j)}(x)=\sum_{l=1}^n\left(L^{(j)}_f(X_l)-\bE_f\left[L^{(j)}_f(X_l)\right]\right)}$,  $j=2,3$, where
$$
L^{(2)}_f(X_l):=\int_{\bR^2}n^{-1}K_{h}\big(p_1\Omega\Gamma y+p_2X_l-\Omega\Gamma QD\Omega x\big)f(y)\rd y,
$$
$$
L^{(3)}_f(X_l):=\int_{\bR^2}n^{-1}K_{h}\big(p_1\Omega\Gamma X_l+p_2y-\Omega\Gamma QD\Omega x\big)f(y)\rd y,
$$
\begin{eqnarray*}
&& |L^{(j)}_f(X_l)|\leq L^2\|\cK\|_1^2(n\delta_n^2)^{-1}\;\;\text{and}
\\*[2mm]
&&\sum_{l=1}^n\text{Var}_f\left[L^{(j)}_f(X_l)\right]\leq 4L^4\left(\|\cK\|_1^4\vee\|\cK\|_2^4\right)n^{-1},\;\text{since}\;p_1^2+p_2^2=1.
\end{eqnarray*}

\noindent Put $\displaystyle{\lambda_p^{(2)}(\epsilon):=2\sqrt{\epsilon^{-1}}\left[\sqrt{2}+\sqrt{3p}\right]\left(\|\cK\|_1^2\vee\|\cK\|_2^2\right)}$.\\ As for $j=2,3$ the $L^{(j)}_f(X_l)$'s are independent variables, we get from Bernstein inequality that for any $p\geq 1$, integer $n\geq3$,  any $z>0$ satisfying $0\leq z\leq 4p\alpha_n$ and any $h>0$ satisfying $h\leq\epsilon^{-1}$ and $\epsilon\alpha_n h\leq n\delta_n^4$
\begin{eqnarray}
\label{uperbound5}
\sup_{j=2,3}\bP_f\bigg\{\sup_{D,Q\in\cQ_{\delta_n},\;Q\neq D}\bigg[\left|\xi_{h,(D,Q)}^{(j)}(x)\right|-\lambda_p^{(2)}(\epsilon)L^2\sqrt{\frac{2p\alpha_n+z}{nh}}\;\bigg]>0\bigg\}
\\
\leq 2 e^{-z},\nonumber
\\
\label{uperbound6}
 \sup_{j=2,3}\bE_f\bigg\{\sup_{D,Q\in\cQ_{\delta_n},\;Q\neq D}\bigg[\left|\xi_{h,(D,Q)}^{(j)}(x)\right|-\lambda_p^{(2)}(\epsilon)L^2\sqrt{\frac{2p\alpha_n+z}{nh}}\;\bigg]_+\bigg\}^p
\\
\leq C_2^{(p)}(\cK,\epsilon)L^{2p}\left[nh\right]^{-\frac{p}{2}} e^{-z},\nonumber
\end{eqnarray}
where $
 C_2^{(p)}(\cK,\epsilon)=2^{\frac{p}{2}+1}\left(3\sqrt{\epsilon^{-1}}\right)^p\Gamma(p+1)\left(\|\cK\|_1^2\vee \|\cK\|_2^2\right)^{p}.$

On the other hand, choosing $\epsilon=1$ and $\alpha_n=\alpha\ln(n)$, one has, in view of (\ref{eq:delta}), for all $p\geq 1$ and all integer $n\geq 3$

\vskip0.1cm

\noindent $\displaystyle{\bP_f\bigg\{\sup_{h\in\mH}\sup_{D,Q\in\cQ_{\delta_n},\;Q\neq D}\left[\left|\xi_{h,(D,Q)}^{(j)}(x)\right|-\alpha\sqrt{6p}\lambda_{p}^{(2)}(1)L^2\sqrt{\ln(n)/n}\;\right]>0\bigg\}}$

\vskip0.1cm

\noindent $\displaystyle{\leq 2\ln(n)n^{-4p}, j=2,3,}$ and then
\begin{eqnarray}
\label{uperbound7}
 \sup_{j=2,3}\bP_f\bigg\{\sup_{h\in\mH}\sup_{D,Q\in\cQ_{\delta_n},\;Q\neq D}\bigg[\left|\xi_{h,(D,Q)}^{(j)}(x)\right|-\frac{1}{3}\widehat{\cU}_n\sqrt{\ln(n)/nh}\;\bigg]>0\bigg\}
\\
\leq 2\ln(n)n^{-4p}+\;\bP_f\bigg\{\alpha\sqrt{6p}\lambda_{p}^{(2)}(1)L^2\sqrt{1/\ln(\ln(n))}>\frac{1}{3}\widehat{\cU}_n\geq \frac{1}{3}\bigg\},\nonumber
\end{eqnarray}
where the second term of the right hand side is equal to zero for $n$ large enough. Similarly we get for all $p\geq 1$ and all integer $n$ large enough
\begin{eqnarray}
\label{uperbound8}
\sup_{j=2,3}\bE_f\bigg\{\sup_{h\in\mH}\sup_{D,Q\in\cQ_{\delta_n},\;Q\neq D}\bigg[\left|\xi_{h,(D,Q)}^{(j)}(x)\right|-\frac{1}{3}\widehat{\cU}_n\sqrt{\ln(n)/nh}\;\bigg]_+\bigg\}^p
\\
 \leq C_2^{(p)}(\cK,1)L^{2p}[\ln(n)]^pn^{-4p}.\nonumber
\end{eqnarray}

\vskip0.1cm

\noindent Now we derive upper bound of $\xi_{h,(D,Q)}^{(1)}(x)$ from exponential inequalities
 developed in \cite{patricia2003}, Theorem $3.4$.

Set  $\;\displaystyle{\xi_{h,(D,Q)}^{(1)}(x)=\sum_{k=2}^n\sum_{l=1}^{k-1}g(X_k,X_l)}$, where
\begin{eqnarray*}
g(X_k,X_l) :=  \varphi(X_k,X_l)+\varphi(X_l,X_k)-\bE_f\left[\varphi(X_k,X_l)+\varphi(X_l,X_k)|X_l\right]
\\
-\;\bE_f\left[\varphi(X_k,X_l)+\varphi(X_l,X_k)|X_k\right]+\bE_f\left[\varphi(X_k,X_l)+\varphi(X_l,X_k)\right].
\end{eqnarray*}
Note that $\displaystyle{\bE_f\left[g(X_k,X_l)|X_l\right]=\bE_f\left[g(X_k,X_l)|X_k\right]=0}$ and
\begin{eqnarray*}
\left|g(X_k,X_l)\right| & \leq & 12(1\vee\|\cK\|_\infty)^2(nh)^{-2}=:A,
\\
\sum_{k=2}^n\sum_{l=1}^{k-1}\bE_f\left[g(X_k,X_l)^2\right] & \leq & 90\left(\|\cK\|_1^4\vee\|\cK\|_2^4\right)L^4(nh)^{-2}=:C^2.
\end{eqnarray*}
Moreover for any $a_k(\cdot), b_k(\cdot), k \in\bN^*,$ verifying  ${\bE_f\big[\sum_{k=2}^na_k(X_k)^2\big]\leq 1}$ and ${\bE_f\big[\sum_{l=1}^{n-1}b_l(X_l)^2\big]\leq 1}$ one has using $2ab\leq a^2+b^2$
\begin{eqnarray*}
&&\hskip-0.65cm\bE_f\bigg[\sum_{k=2}^n\sum_{l=1}^{k-1}g(X_k,X_l)a_k(X_k)b_l(X_l)\bigg]\leq 4(n-1)\sup_{u\in\bR^2}\bE_f\left|\varphi(u,X_1)+\varphi(X_1,u)\right|
\\
&&\hskip5.41cm\leq8L^2\|\cK\|_1^2(n\delta_n^2)^{-1}=:D.
\end{eqnarray*}
 By independence of the $X_k$'s one has for any $u\in\bR^2$
\begin{eqnarray*}
&&\sum_{l=1}^{k-1}\bE_f\left[g(u,X_l)^2|X_k\right]
\\
&& = \sum_{l=1}^{k-1}\text{Var}_f\left[\left(\varphi(u,X_l)+\varphi(X_l,u)\right)-(n-1)^{-1}\left(L^{(2)}_f(X_l)+L^{(3)}_f(X_l)\right)\right]
\\
&& \leq 24L^4\left(\|\cK\|_1^4\vee\|\cK\|_2^4\right)/n^3h^2\delta_n^2.
\end{eqnarray*}
Similarly, one has

\vskip0.1cm

\centerline{$\displaystyle{\sup_{l=1,\ldots,n-1}\sup_{u\in\bR^2}
\bigg\{\sum_{k=l+1}^{n}\bE_f\left[g(X_k,u)^2|X_l\right]\bigg\}\leq\frac{24L^4\left(\|\cK\|_1^4\vee\|\cK\|_2^4\right)}{n^3h^2\delta_n^2}=:B^2}$.}

\vskip0.1cm

\noindent It gives for any integer $n\geq 3$ and any real number $z>0$

\vskip0.1cm

\centerline{$\displaystyle{\bP_f\left\{\left|\xi_{h,(D,Q)}^{(1)}(x)\right|\geq\cU(z)\right\}\leq 6 e^{-z}}$,}

\vskip0.1cm

\noindent where $\cU(z):=4\sqrt{2}C\sqrt{z}+8\sqrt{2}Dz+426Bz^{3/2}+414Az^2$.

\vskip0.1cm

\noindent By integration of the latter inequality we obtain for all $p\geq 1$, all integers $n\geq 3$ and any $z\geq 1$

\vskip0.1cm

\centerline{$\displaystyle{\bE_f\bigg\{\left[\left|\xi_{h,(D,Q)}^{(1)}(x)\right|-\cU(z)\right]_+\bigg\}^p\leq 3\times 2^{2p+1}\Gamma(2p+1)\left[z\cU(1)\right]^p e^{-z}}$.}

\vskip0.1cm

\noindent Put $\displaystyle{\lambda_p^{(3)}(\epsilon):=87955p\sqrt{p}\big(\|\cK\|_1^2\vee\|\cK\|_2^2\vee\|\cK\|_\infty^2\big)\Big[\sqrt{\epsilon^{-1}}\vee\epsilon^{-1}\vee\epsilon^{-\frac{3}{2}}\Big]}$.\\ It follows that for all $p\geq 1$, all integer $n\geq 3$, all real number $z$ satisfying $0<z\leq 4p\alpha_n$ and all real number $h>0$ satisfying $nh\geq\epsilon\alpha_n$, $\epsilon\alpha_n h\leq n\delta_n^4$ and $\epsilon\alpha_n\leq n\delta_n\sqrt{h}$
\begin{gather}
\label{uperbound9}
\bP_f\bigg\{\sup_{D,Q\in\cQ_{\delta_n},\;Q\neq D}\bigg[\left|\xi_{h,(D,Q)}^{(1)}(x)\right|-\lambda_p^{(3)}(\epsilon)L^2\sqrt{\frac{2p\alpha_n+z}{nh}}\;\bigg]>0\bigg\}
\\
\leq 6 e^{-z},
\nonumber
\\
\label{uperbound10}
\bE_f\bigg\{\sup_{D,Q\in\cQ_{\delta_n},\;Q\neq D}\left[\left|\xi_{h,(D,Q)}^{(1)}(x)\right|-\lambda_p^{(3)}(\epsilon)L^2\sqrt{\frac{2p\alpha_n+z}{nh}}\;\right]_+\bigg\}^p
\\
\leq C_3^{(p)}(\cK,\epsilon)6\times 2^{2p}\Gamma(2p+1)L^{2p}\left[nh/\alpha_n\right]^{-\frac{p}{2}} e^{-z},
\\*[2mm]
C_3^{(p)}(\cK,\epsilon)=
\left(\|\cK\|_1^2\vee\|\cK\|_2^2\vee\|\cK\|_\infty^2\right)^p\left[45038p\left(\sqrt{\epsilon^{-1}}\vee\epsilon^{-1}\vee\epsilon^{-\frac{3}{2}}\right)\right]^{p}.\nonumber
\nonumber
\end{gather}


\noindent In another hand, as previously, we get for all $p\geq 1$ and all $n$ large enough
\begin{eqnarray}
\label{uperbound11}
\bP_f\bigg\{\sup_{h\in\mH}\sup_{D,Q\in\cQ_{\delta_n},\;Q\neq D}\bigg[\left|\xi_{h,(D,Q)}^{(1)}(x)\right|-\frac{1}{3}\widehat{\cU}_n\sqrt{\ln(n)/nh}\;\bigg]>0\bigg\}
\\
\leq 6\ln(n)n^{-4p},\nonumber
\end{eqnarray}
\begin{eqnarray}
\label{uperbound12}
\bE_f\bigg\{\sup_{h\in\mH}\sup_{D,Q\in\cQ_{\delta_n},\;Q\neq D}\bigg[\left|\xi_{h,(D,Q)}^{(1)}(x)\right|-\frac{1}{3}\widehat{\cU}_n\sqrt{\ln(n)/nh}\;\bigg]_+\bigg\}^p
\\
\leq C_3^{(p)}(\cK,1)L^{2p}[\ln(n)]^pn^{-4p}.\nonumber
\end{eqnarray}

\smallskip

\noindent\textit{Third step: end of proofs of Propositions \ref{prop1}-\ref{prop4}.}

\vskip0.1cm

\noindent Remind that third assertion of Proposition \ref{prop1} is already proved in step one. First and second ones follow from inequalities (\ref{uperbound4}), (\ref{uperbound7}), (\ref{uperbound11}) and (\ref{uperbound3}), (\ref{uperbound8}), (\ref{uperbound12}) respectively, since $\boldsymbol{\ma}\geq 1$.

\vskip0.1cm

\noindent End of the proof of Proposition \ref{prop4}. Note first that $\mathbf{B}\geq \mathbf{B}_1\vee\mathbf{B}_2\vee\mathbf{B}_3$, where
\begin{eqnarray*}
\mathbf{B}_1&=&20p\big[1+\epsilon_1^{-1/2}\big]\big(\lambda_{2p}^{(1)}(\epsilon_1)\|\cK\|_1\big)^2L^2+8C(\cK,\mb,\sqrt{2})L^2;
\\
\mathbf{B}_2&=&6\sqrt{6p}\lambda_{p}^{(2)}(\epsilon_2)L^2+8C(\cK,\mb,\sqrt{2})L^2;
\\
\mathbf{B}_3&=&6\sqrt{6p}\lambda_{p}^{(3)}(\epsilon_3)L^2+8C(\cK,\mb,\sqrt{2})L^2
\end{eqnarray*}
Note first that for any $\mfi=1,\ldots,\mfi^*$

\vskip0.1cm

\centerline{$\displaystyle{L^2h_\mfi^\beta=\sqrt{\frac{\omega_\mfi}{\bn h_\mfi}},\;\;\bn h_\mfi\geq\left[\left(8+4\alpha\right)L^{\frac{2}{\beta}}\right]^{-\frac{2\beta}{2\beta+1}}\omega_\mfi,\;\;\;
\ln\left(\frac{\omega_{\mfi-1}\bn}{\mathbf{n}_{\mfi-1}\omega_\mfi}\right)\leq 3\omega_\mfi}$.}

\vskip0.1cm

\noindent Thus, in view of (\ref{uperbound1}) and (\ref{uperbound2}) with $\displaystyle{\epsilon=\epsilon_1=\big[\left(8+4\alpha\right)L^{\frac{2}{\beta}}\big]^{-\frac{2\beta}{2\beta+1}}}$ and $\alpha_n=\omega_\mfi$ one has for all $p\geq 1$ and all integer $n\geq 3$
\begin{eqnarray*}
&&\sup_{b\in\{d,d_\perp\}}\bP_f\bigg\{\sup_{D\in\cQ_{\delta_n}}\left[\left|\xi^{(\mfi)}_{(h_\mfi,b)}(x)\right|-
\lambda_{2p}^{(1)}(\epsilon_1)\sqrt{10p}\|\cK\|_1L^3h_\mfi^\beta\;\right]>0\bigg\}
\\
&&\leq \frac{2}{e^{8}}\left(\frac{\omega_{\mfi-1}\bn}{\mathbf{n}_{\mfi-1}\omega_\mfi}\right)^{\frac{p\beta}{2\beta+1}},
\\
&&\sup_{b\in\{d,d_\perp\}}\bE_f\bigg\{\sup_{D\in\cQ_{\delta_n}}\left[\left|\xi^{(\mfi)}_{(h_\mfi,b)}(x)\right|-
\lambda_p^{(1)}(\epsilon_1)\sqrt{2p}\|\cK\|_1L^3h_\mfi^\beta\;\right]_+\bigg\}^p
\\
&&\leq C_1^{(p)}(\cK,\epsilon_1)L^{3p}h_\mfi^{p\beta}\;\;\text{and}
\\*[1mm]
&&\sup_{b\in\{d,d_\perp\}}\bE_f\Big\{\sup_{D\in\cQ_{\delta_n}}\left|\xi^{(\mfi)}_{(h_\mfi,b)}(x)\right|\;\Big\}^p
\\
&&\leq 2^{p-1}\left[C_1^{(p)}(\cK,\epsilon_1)+\left(\lambda_p^{(1)}(\epsilon_1)\sqrt{2p}\|\cK\|_1\right)^p\right]L^{3p}h_\mfi^{p\beta}.
\end{eqnarray*}
By  Cauchy-Schwartz inequality, noting that
$$
|\xi^{(\mfi)}_{h_\mfi,D}(x)|\leq\sqrt{2}\|\cK\|_1L\left(|\xi^{(\mfi)}_{(h_\mfi,d)}(x)|+|\xi^{(\mfi)}_{(h_\mfi,d_\perp)}(x)|\right)+|\xi^{(\mfi)}_{(h_\mfi,d)}(x)|\times|\xi^{(\mfi)}_{(h_\mfi,d_\perp)}(x)|
$$
we easily get for all $p\geq 1$, all $\mfi=1,\ldots\mfi^*$ and all integer $n\geq 3$
\begin{eqnarray}
\label{uperbound13}
\bP_f\bigg\{\sup_{D\in\cQ_{\delta_n}}\left[|\xi^{(\mfi)}_{h_\mfi,D}(x)|-\mathbf{C} L^2h_\mfi^\beta\;\right]>0\bigg\}\leq \frac{8}{e^{8}}\left(\frac{\omega_{\mfi-1}\bn}{\mathbf{n}_{\mfi-1}\omega_\mfi}\right)^{\frac{p\beta}{2\beta+1}}\;\text{and}
\\
\label{uperbound14}
\bE_f\left\{\sup_{D\in\cQ_{\delta_n}}\left[|\xi^{(\mfi)}_{h_\mfi,D}(x)|-\mathbf{C} L^2h_\mfi^\beta\;\right]_+\right\}^p\leq C'h_\mfi^{p\beta},
\end{eqnarray}
where$C':=C'(p,\cK,\beta,L,\alpha)>0$ and $\mathbf{C}=2^{-1}\mathbf{B}-4C\big(\cK,\mb,\sqrt{2}\big)L^2$.

\noindent Remark now that for any $\mfi=1,\ldots,\mfi^*$ one has $\displaystyle{h_\mfi\leq\left[L^{-4}(8+4\alpha)\right]^{\frac{1}{2\beta+1}}}$ and

\vskip0.1cm

\centerline{$\displaystyle{\omega_\mfi h_\mfi\leq L^{-\frac{4}{2\beta+1}}(9+3\alpha)^{\frac{2\beta+2}{2\beta+1}}C(\beta)\bn\delta_n^4}$.}

\vskip0.1cm

\noindent Thus, in view of (\ref{uperbound5}) and (\ref{uperbound6}) with $\displaystyle{\epsilon=\epsilon_2=L^{\frac{4}{2\beta+1}}(9+4\alpha)^{-\frac{2\beta+2}{2\beta+1}}\left[C(\beta)\right]^{-1}}$ and $\alpha_n=\omega_\mfi$ one has for all $p\geq 1$, all $\mfi=1,\ldots,\mfi^*$ and all integer $n\geq 3$
\begin{eqnarray}
\label{uperbound15}
\sup_{j=2,3}\bP_f\bigg\{\sup_{D,Q\in\cQ_{\delta_n},\;D\neq Q}\left[|\xi^{(j),(\mfi)}_{h_\mfi,D}(x)|-\frac{1}{3}\mathbf{C} L^2h_\mfi^\beta\;\right]>0\bigg\}
\\
\leq \frac{2}{e^{8}}\left(\frac{\omega_{\mfi-1}\bn}{\mathbf{n}_{\mfi-1}\omega_\mfi}\right)^{\frac{p\beta}{2\beta+1}};\nonumber
\\
\label{uperbound16}
\sup_{j=2,3}\bE_f\bigg\{\sup_{D,Q\in\cQ_{\delta_n},\;D\neq Q}\left[|\xi^{(j),(\mfi)}_{h_\mfi,D}(x)|-\frac{1}{3}\mathbf{C} L^2h_\mfi^\beta\;\right]_+\bigg\}^p\leq C''h_\mfi^{p\beta},
\end{eqnarray}
where $C'':=C''(p,\cK,\beta,L,\alpha)>0$.
 Note finally that for any $\mfi=1,\ldots,\mfi^*$

\vskip0.1cm

\centerline{$\displaystyle{\omega_\mfi\leq L^{\frac{2}{2\beta+1}}(9+3\alpha)^{\frac{2\beta+1/2}{2\beta+1}}C(\beta)\bn\delta_n\sqrt{h_\mfi}}$.}

\vskip0.1cm

\noindent Thus, in view of (\ref{uperbound9}) and (\ref{uperbound10}) with $\displaystyle{\epsilon=\epsilon_3=L^{-\frac{4}{2\beta+1}}(9+4\alpha)^{-\frac{2\beta+2}{2\beta+1}}\left[C(\beta)\right]^{-1}}$ and $\alpha_n=\omega_\mfi$ one has for all $p\geq 1$, all $\mfi=1,\ldots,\mfi^*$ and all integer $n\geq 3$
\begin{eqnarray}
\label{uperbound15}
&&\bP_f\bigg\{\sup_{D,Q\in\cQ_{\delta_n},\;D\neq Q}\left[|\xi^{(1),(\mfi)}_{h_\mfi,D}(x)|-\frac{1}{3}\mathbf{C} L^2h_\mfi^\beta\;\right]>0\bigg\}
\\
&&\leq 6e^{-8}\bigg(\frac{\omega_{\mfi-1}\bn}{\mathbf{n}_{\mfi-1}\omega_\mfi}\bigg)^{\frac{p\beta}{2\beta+1}};\nonumber
\\
\label{uperbound16}
&&\bE_f\bigg\{\sup_{D,Q\in\cQ_{\delta_n},\;D\neq Q}\left[|\xi^{(1),(\mfi)}_{h_\mfi,D}(x)|-\frac{1}{3}\mathbf{C} L^2h_\mfi^\beta\;\right]_+\bigg\}^p\leq C'''h_\mfi^{p\beta},
\end{eqnarray}
where $C''':=C'''(p,\cK,\beta,L,\alpha)>0$.  Proposition \ref{prop4} is proved.
\epr

\bibliographystyle{agsm}

\end{document}